\documentclass[a4paper,11pt]{article}
\usepackage{cancel}
\usepackage{amsmath}
\usepackage{amssymb}
\usepackage{soul}
\usepackage{arydshln}
\usepackage{adjustbox, lipsum}
\usepackage{enumerate}
\usepackage[ruled,vlined,linesnumbered]{algorithm2e}
\usepackage{rotating}
\usepackage[utf8]{inputenc}
\usepackage{amsfonts}
\usepackage{pstricks}
\usepackage{pdflscape}
\usepackage{graphicx}
\usepackage{soul}
\usepackage{afterpage}
\usepackage{rotating}
\usepackage{natbib}
\usepackage{booktabs}
\usepackage{authblk}
\usepackage{multirow}
\usepackage{mathtools}
\usepackage{scalerel}

\newtheorem{teor}{\bf Theorem}[section]

\newtheorem{lema}{\bf Lemma}[section]

\newtheorem{prop}{\bf Proposition}[section]
\newtheorem{coro}{\bf Corollary}[section]
\newtheorem{ejem}{\bf Example}[section]

\SetKwRepeat{Do}{do}{while}
\providecommand{\norm}[1]{\lVert#1\rVert}
\DeclarePairedDelimiter\abs{\lvert}{\rvert}%

\newcommand{\dem}{\par \noindent{\bf Proof:} }
\newcommand{\fin}{\hfill $\square$  \par \bigskip}
\newcommand{\bb}{b}
\usepackage{xcolor}
\definecolor{gr}{RGB}{0,153,0}

\author[(a)]{Marta Baldomero-Naranjo\footnote{Corresponding author}}
\author[(a)]{Luisa I. Martínez-Merino}
\author[(a)]{Antonio M. Rodríguez-Chía}

\affil[(a)]{\small{Departamento de Estadística e Investigación Operativa, Universidad de Cádiz, Cádiz, Spain, marta.baldomero@uca.es, luisa.martinez@uca.es, antonio.rodriguezchia@uca.es}}

\date{}

\title{\LARGE Tightening big Ms in Integer Programming Formulations for Support Vector Machines with Ramp Loss}

\setcitestyle{authoryear,open={(},close={)}}

\makeatletter
\let\oldabs\abs
\def\abs{\@ifstar{\oldabs}{\oldabs*}}
\makeatother
\begin{document}

%
\maketitle
\begin{abstract}
	 This paper considers various models of support vector machines with ramp loss, these being an efficient and robust tool in supervised classification for the detection of outliers. The exact solution approaches for the resulting optimization problem are of high demand for large datasets. Hence, the goal of this paper is to develop algorithms that provide  efficient methodologies to exactly solve these optimization problems. These approaches are based on three strategies for obtaining tightened values of the big M parameters included in the formulation of the problem. Two of them require solving a sequence of continuous problems, while the third uses the Lagrangian relaxation to tighten the bounds. The  proposed resolution methods are valid for the $\ell_1$-norm and $\ell_2$-norm ramp loss formulations. They were tested and compared with existing solution methods 
	in simulated and real-life datasets, showing the efficiency of the developed methodology. 

\noindent		
	\textbf{Keywords:} Location; Support Vector Machine; Ramp Loss Model; Mixed Integer Programming; Indicator Constraints. 

\end{abstract}

\section{Introduction}
\label{sec:introduction}

The goal of supervised  classification is to provide a decision rule to classify  the new individuals of a population in classes,
based on the learning over a set of individuals already classified in those classes. The field of  mathematical programming has resulted in a very useful tool for supervised classification, introduced by \cite{Mangasarian65,Mangasarian68}. In particular, one of its most useful approaches are the support vector machines (SVM). Since their introduction by \cite{13SVN} and \cite{Vapnikunicidad}, SVM have been applied in many fields, including: biology \citep{ap9},  
medicine \citep{ap8,ap4,apCancer}, bioinformatics \citep{ap1,ap10}, machine vision \citep{ap6,ap2}, text classification \citep{ap3},  
and finance \citep{ap7,apcreditscoring,ap5}, 
among others. 

In order to define a supervised classification problem more precisely, we consider the following notation. Given a set $N$ of individuals partitioned into two classes $\cal{Y}=\mbox{\{-1,1\}}$, each individual $i \in N$ is associated with a pair $(x_i,y_i) \in \mathbb{R}^d \times \{-1,1\}$, where $d$ is the number of features analyzed in each individual of $N,$ $x_i$ contains the feature values, and $y_i$ provides the class membership (1 or -1). The aim of the SVM is to build a separating hyperplane $w\cdot x + \bb=0 $ to classify new individuals. The search for this hyperplane is based on a compromise between maximizing the distance (margin) between the two parallel hyperplanes supporting individuals of each class and minimizing the error caused by misclassified data.
As proven in \cite{03carrioza2013}, the maximization of the margin between the parallel hyperplanes associated with each class can equivalently be expressed as the minimization of $\norm{w}^\circ,$ (where $\norm{\cdot}^\circ$ is the dual norm of the norm used to measure the distance to the hyperplanes in $ \mathbb{R}^d$). In \cite{BlancoPuertoChia}, the methodology applied to classical SVM (using $\ell_2$-norm) is extended to the general $\ell_p$-norm with $p\geq 1$. Although we will concentrate on cases $p=1$ and $p=2$, we provide the following common formulation for both cases that it is also valid for any value $p\geq 1$,
\begin{alignat}{4}
\mbox{(SVM-$\ell_p$)} & \quad &  & \mbox{min}&  \quad & 	\displaystyle \frac{1}{p}\left(\norm{w}_p\right)^ p + C  \sum_{i=1}^n\xi_i , & & \nonumber \\
& & &\mbox{s.t.} &  \quad &  y_i \left(\sum_{k=1}^{d}w_k x_{ik}+\bb \right) \geq 1- \xi_i,   &  & i \in N,  \\
& & & & &  \xi_i \geq 0, & \quad  & i \in N,
\end{alignat}
where the $w$- and $\bb$-variables are the separation hyperplane coefficients and the $\xi$-variables represent classification errors. $C$ is a non-negative constant, $N$ is the set $\{1,\dots,n\}$,  and $\norm{w}_p$ represents the $\ell_p$-norm of vector $w$  with $p\geq 1$.  Therefore, the objective function given above seeks to maximize the margin and minimize the sum of misclassified data deviations, given by $\sum_{i=1}^n\xi_i$. Parameter $C$ regulates the trade-off between both goals. 
  
As shown, the classical formulation of SVM uses an unbounded continuous variable for measuring the deviations of misclassified data: $\xi_i,$ for $i \in N$.  This results in models that are sensitive to outliers. In order to avoid that, a more robust classifier  against outliers has been proposed in the literature: SVM with ramp loss. This was introduced by \cite{31Brooks2011} for a mixed integer problem with conditional constraints. In this model, the deviation of misclassified data is truncated to avoid extreme values: i.e., for a given deviation, the penalization is bounded. Due to its usefulness as a classifier, this model has attracted the attention of many researchers and  has been analyzed whilst considering the $\ell_1$-norm and $\ell_2$-norm (see \cite{31Brooks2011, carrizosaheuristic,bonami,Blanco2018}). This model has also been studied under the framework of statistical learning theory (see \cite{RLLPSVM}).

  In spite of the interest of this model for data  classification avoiding the effect of outliers, the search for exact solution approaches for large datasets is still an ongoing problem, as stated in \cite{PedroDuarte17}. 
In this work, we focus our attention on this challenge by providing efficient algorithms to obtain exact optimal solutions for the resulting optimization problems from the ramp loss model.
We first consider the model using the $\ell_1$-norm 
 and then the presented strategies are adapted to the $\ell_2$-norm. 
The analysis of these two norms is fully justified because the $\ell_1$-norm has various interesting properties, such as its sparsity, whilst the $\ell_2$-norm is the one used in classic SVMs. Indeed, using the $\ell_1$-norm when the dataset contains many features is  worthwhile as it results in easily interpreted classifiers with a reduced number of selected features (see \cite{23Antonio,22Luisa,21Maldonado2014}). 

\cite{31Brooks2011} introduces a new integer programming formulation for SVM with ramp loss that accommodate the use of nonlinear kernel functions. However, he proved that SVM with ramp loss can produce robust classifiers when using a linear  kernel in the presence of outliers.  For this reason, we concentrate on SVM ramp loss models using a linear kernel.

 The developed methodologies in this paper are compared with the methods proposed in \cite{bonami}.
 In particular, \cite{bonami} introduce a non convex formulation for the SVM with ramp loss using the $\ell_2$-norm. They also propose two methods based on bound tightening approaches for efficiently solving the model. One of the methods is focused on an iterative tightening of the bounds of variables providing locally valid bounds for the big M parameters appearing in the model. This method is included in {CPLEX} arsenal and it is {known} as {\sl local implied bound cuts}. The second method tightens  the bounds of the variables iteratively by solving MIPs.  
It should be remarked that the model studied in \cite{bonami} provide an unique initial big M parameter which is valid for all the family of constraints in which these parameters appear. 

In contrast, this work proposes strategies that provide different tightened bounds on each big M parameter. In addition, unlike \cite{bonami}, our strategies are based on the tightening of bounds by iteratively solving linear programming problems in the $\ell_1$-norm case and quadratic ones in the $\ell_2$-norm, but non integer. In the computational results section we will analyze the positive effect of our approaches. 

 The remainder of the paper is organized as follows. Section \ref{general_model} introduces the model and presents a set of valid inequalities for the formulation. 
In Section \ref{sec:Strategies}, 
   tightened values for the big M parameters in the model with $\ell_1$-norm are proven and strategies for obtaining them are proposed. Some of these strategies are based on tightening values of the $w$-variables and others are based on the Lagrangian relaxation of the model. In Section \ref{sec:l2norm}, 
   all the previously proposed strategies are adapted to the $\ell_2$-norm model. Section \ref{sec:computationalExp}
    contains computational experiments carried out on simulated and real-life datasets. Our conclusions and other potential research topics are included in Section \ref{conclusions}.  Lastly, the proofs of the results for the $\ell_2$-norm model can be found in the Appendix \ref{sec:Appendix}.
    
\section{The model}
\label{general_model}

We will focus our research on the Support Vector Machine model with ramp loss, as introduced by \cite{31Brooks2011}. The model which uses the $\ell_p$-norm for any $p\geq 1$ is formulated as a mixed integer  program with conditional constraints:

\begin{alignat}{4}
\mbox{(RL-$\ell_p$)} &\quad  &  & \mbox{min}&  \quad & 	\displaystyle \frac{1}{p}(||w||_p)^p + C  \left( \sum_{i=1}^n\xi_i + 2 \sum_{i=1}^n z_i \right), & & \nonumber \\
& & &\mbox{s.t.} &  \quad &  \mbox{if } z_i=0,\;\; y_i \left(\sum_{k=1}^{d}w_k  x_{ik}+\bb \right) \geq 1- \xi_i,   &  & i \in N, \label{rl-cons} \\
& & & & & 0 \leq   \xi_i \leq 2, & \quad  & i \in N, \label{rl-l2-ec2}  \\
& & & & & z_i\in \{0,1\}, & &  i \in N, \label{rl-l2-ec3} 
\end{alignat}

In this model, $\xi_i$ determines the penalization if the misclassified object $i$ is  within the strip defined by the two parallel hyperplanes ($wx+b=1$ and $wx+b=-1$) and $z_i$ determines whether $i$ is a misclassified object outside this strip or not.
This model can be linearized by replacing the set of conditional constraints \eqref{rl-cons} with big M constraints as follows,
\begin{alignat}{4}
\mbox{(RL-$\ell_p$-M)} & \quad &  & \mbox{min}&  \quad & 	\displaystyle \frac{1}{p}(||w||_p)^p + C  \left( \sum_{i=1}^n\xi_i + 2 \sum_{i=1}^n z_i \right), & & \nonumber \\
& & &\mbox{s.t.} &  \quad & \eqref{rl-l2-ec2},\eqref{rl-l2-ec3},   & \quad & \nonumber \\
& & & & & y_i\left( \sum_{k=1}^{d}w_k x_{ik}+\bb\right) \geq 1- \xi_i -M_iz_i,   & \quad  & i \in  N, \label{ec:linealizado_lp}
\end{alignat}
 where $M_i$ is a big enough constant, for $i \in N$. In the following, we give a result that establishes a relationship between the $\xi$-variables and the $z$-variables.

\begin{prop}\label{remark1}
	
	\begin{enumerate}[i)]
		\item An optimal solution of (RL-$\ell_p$), with $p\geq 1,$ $(w^*, \bb^*, \xi^*, z^*),$ satisfies the following condition:
		\begin{equation}
		\xi_i^* z_i^*=0, \label{cond:1} \quad i \in N. 
		\end{equation}
	
		\item Using $M_i,$ for $i \in N,$  so that  an optimal solution of (RL-$\ell_p$-M) satisfies \eqref{cond:1} does not imply that (RL-$\ell_p$-M) and (RL-$\ell_p$) are equivalent. 
	\end{enumerate}
	
\end{prop}
The proof of \textit{i)} can be obtained using a simple contradiction, and for this reason it has been omitted.  An example showing  \textit{ii)} is given in the following for $p=1$ and $p=2$.

	\begin{ejem}
Let  $N$ be a set of individuals partitioned into two classes, where the associated pair $(x_i,y_i)$ for $i\in  N$ is $\{((-2,1),1),((-1,-1),1),((-5,-3),-1),((1,3),-1),((1,0),-1)\}$ and $C=10.$  The optimal solution of (RL-$\ell_1$-M) is 23, a valid value for M is 41, and  the optimal solution is ${w}^*=(-3,0)$, $b^*=2,$ $\xi^* =(0,0,0,0,0),$ and $z^*=(0,0,1,0,0),$ (the valid value for the big M parameter was computed following the strategies described in Section \ref{sec:Strategies}).  Nevertheless, considering $M=5$ in model (RL-$\ell_1$-M), the optimal value is 25.8333; where ${w}^*=(-2.5,0)$, $b^*=1.5,$ $\xi^* =(0,0.3333,0,0,0),$ and $z^*=(0,0,1,0,0);$ i.e., condition \eqref{cond:1} is still fulfilled for  $i \in N$, but its optimal objective value differs. 
Similarly, the optimal solution of (RL-$\ell_2$-M) is 23.6, a valid value for M is 14, and the optimal solution is $w^*=(-2.4, -1.2), \bb^*=1.4, \xi=(0,0,0,0,0),$ and $z=(0,0,1,0,0),$ ( the valid value for the big M parameter was computed following the strategies described in Section \ref{sec:l2norm}). However, establishing $M=5$ in model (RL-$\ell_2$-M), the optimal value is 25.9333; where $w^*=(-2.2, -0.6), \bb^*=1.2,  \xi^*=(0,0.3333,0,0,0),$ and $z=(0,0,1,0,0).$
\end{ejem}

As a consequence of the previous result, the linearized version of condition \eqref{cond:1} is given by 
	\begin{equation}\label{rest:cond1}
	\xi_i \leq 2(1 - z_i), \quad i \in N.
	\end{equation}
 This set of constraints will be used to strengthen the (RL-$\ell_p$-M) formulation. Henceforth, we will refer to \mbox{(RL-$\ell_p$-M)}$+\eqref{rest:cond1}$ as (RL-$\ell_p$-M), unless it is stated otherwise.

Observe that in the (RL-$\ell_p$-M) model,  the choice of an appropriate value for $M_i,$ for $i \in N$ is essential when providing efficient solution approaches. Note that $M_i$ should be big enough for the equivalence between (RL-$\ell_p$) and (RL-$\ell_p$-M),  but it should also be as small as possible so as to improve the linear relaxation and the computational time needed to solve it. 
 The following proposition provides valid values for the big M parameters considering a general $\ell_p$-norm, with $p\geq 1$.

\begin{prop}\label{prop1}
For a given $p\geq 1$, the problems (RL-$\ell_p$) and (RL-$\ell_p$-M) are equivalent if 
	$$M_i \geq \left(\underset{j \in N}{\max } \; \{ \norm{x_i-x_j}_{\bar{q}} : y_i=y_j \} \right) \norm{w}_{\bar{p}}, \;\;\mbox{for } i \in N,$$ 
 with $\bar{p},\bar{q}\geq 1$ such that $\frac{1}{\bar{p}}+\frac{1}{\bar{q}}=1$ (i.e. $\norm{\cdot}_{\bar{p}}$ and $\norm{\cdot}_{\bar{q}}$ are dual norms).
\end{prop}

\dem
Taking into account set of constraints \eqref{ec:linealizado_lp}, it holds that: 
$$M_iz_i \geq -y_i\left( \sum_{k=1}^{d}w_k x_{ik}+\bb\right) + 1 - \xi_i, \;\; \mbox{for } i \in N.$$ 
According to Proposition \ref{remark1}, an optimal solution of (RL-$\ell_p$) satisfies condition \eqref{cond:1}. Consequently, a valid value for  $M_i$ in the formulation (RL-$\ell_p$-M) would be one satisfying: 
$$M_i \geq -y_i\left( \sum_{k=1}^{d}w_k x_{ik}+\bb\right) + 1, \;\; \mbox{for } i \in N.$$ 
Hence, a value of $M_i$ satisfying: 
\begin{equation}\label{ec:1}M_i \geq \abs{-y_i\left( \sum_{k=1}^{d}w_k x_{ik}+\bb\right) + 1}=\abs{y_i\left( \sum_{k=1}^{d}w_k x_{ik}+\bb\right) - 1},\end{equation}
would also be valid. On the other hand and using the $\ell_{\bar{q}}$-norm for $\bar{q}\geq 1$, the distance between $x_i$ and the hyperplane $H(y_iw, y_i\bb -1):=\{x: y_i(w\cdot x+\bb) -1=0\},$ for $i \in N,$ (see \cite{Plastria2001}) is: 	
\begin{equation}\label{ec:2}d_{\ell_{\bar{q}}}(x_i, H(y_iw, y_i\bb -1))= \dfrac{\abs{y_i(w\cdot x_i+\bb) -1}}{\norm{w}_{\bar{p}}}.\end{equation}
Thus, taking into account expressions \eqref{ec:1} and \eqref{ec:2}, the (RL-$\ell_p$) and (RL-$\ell_p$-M) problems will be equivalent if the following inequality holds:
\begin{equation}M_i \geq d_{\ell_{\bar{q}}}(x_i, H(y_iw, y_i\bb -1)) \norm{w}_{\bar{p}}, \;\; \mbox{for } i \in N.\end{equation}
Observe that the previous expression represents the distance using the $\ell_{\bar{q}}$-norm from $x_i$ to  $H(y_iw, y_i\bb -1),$ i.e., the supporting hyperplane for each class. Consequently,   $d_{\ell_{\bar{q}}}(x_i, H(y_iw, y_i\bb -1))$ is, at most, the maximum distance between two individuals of the same class, i.e.,
$$d_{\ell_{\bar{q}}}(x_i,H(y_iw, y_i\bb -1))\leq \underset{j \in N}{\max } \;\{ \norm{x_i-x_j}_{\bar{q}} : y_i=y_j \} .$$
Thus, if $M_i$ satisfies	
$$M_i \geq \left(\underset{j \in N}{\max } \; \{ \norm{x_i-x_j}_{\bar{q}} : y_i=y_j \} \right) \norm{w}_{\bar{p}} , \;\; \mbox{for } i \in N,$$ 
both problems will be equivalent. 
\fin

 Strategies for obtaining tighter values than the ones provided by the previous expressions will be presented in the next sections.  Particularly, in Section \ref{sec:Strategies} we will consider the ramp loss model with $\ell_1$-norm while in Section \ref{sec:l2norm}, the attention will be focused on the $\ell_2$-norm. 

\section{Strategies for the $\ell_1$-norm case}
\label{sec:Strategies}

 The objective of this section is to present the model for the $\ell_1$-norm case and to improve the values of big M parameters appearing in the model. As a result, two algorithms are derived. They are based on tightening bounds of $w$-variables in the model and using these bounds to provide tighter values of big M parameter.

A formulation of the SVM with ramp loss using the $\ell_1$-norm is obtained by decomposing the unrestricted variables $w_k$ as the difference of two non-negative variables $w_k^+$ and $w_k^-$ for $k \in D$,  where $D$ is the set $\{1,\dots,d\}$, (see \cite{22Luisa}). In this reformulation, $w_k= w_k^+-w_k^-$, where $ w_k^+, w_k^- \geq 0,$ for $k \in D.$ Thusly,  $\left|w_k\right|= w_k^+ + w_k^-$ in any optimal solution, since $ w_k^+ + w_k^-,$ for  $k \in D$ is part of the objective function to be minimized. This means that, at most, only one of the two variables for any $k\in D$  is non-zero in an optimal solution. The result is the following formulation:
\begin{alignat}{4}
\mbox{(RL-$\ell_1$)} & \quad &  & \mbox{min}&  \quad & 	\displaystyle \sum_{k=1}^{d}(w_k^+ + w_k^-) + C  \left( \sum_{i=1}^n\xi_i + 2 \sum_{i=1}^n z_i \right), & & \nonumber \\
& & &\mbox{s.t.} & \quad & \eqref{rl-l2-ec2}, \eqref{rl-l2-ec3}, \nonumber \\
& & & & &  \mbox{if } z_i=0,\;\; y_i \left(\sum_{k=1}^{d}(w_k^+-w_k^-) x_{ik}+\bb \right) \geq 1- \xi_i,     & \quad &i \in N, \label{ec:cond}  \\
&&& & &  w_k^+\geq 0, w_k^- \geq 0, & \quad &k \in D.   \label{ec:w} 
\end{alignat}
This model can be linearized by replacing set of conditional constraints \eqref{ec:cond} with big M constraints, where $M_i$ is a big enough constant, for $i \in N$. 
\begin{alignat}{4}
\mbox{(RL-$\ell_1$-M)} & \quad &  & \mbox{min}&  \quad & 	\displaystyle \sum_{k=1}^{d}(w_k^+ + w_k^-) + C  \left( \sum_{i=1}^n\xi_i + 2 \sum_{i=1}^n z_i \right), & & \nonumber \\
& & &\mbox{s.t.} &  \quad & \eqref{rl-l2-ec2},\eqref{rl-l2-ec3}, \eqref{ec:w},  & \quad & \nonumber \\
& & & & & y_i\left( \sum_{k=1}^{d}(w_k^+ - w_k^-)x_{ik}+\bb\right) \geq 1- \xi_i -M_iz_i,   & \quad  & i \in N. \label{ec:linealizado}
\end{alignat}

\subsection{Tightening bounds of $\pmb{w}$-variables}

In the previous section, Proposition \ref{prop1} provided valid values $M_i$, for $i\in N,$ depending on $\norm{w}_{\bar p}$.
 Below, we develop two strategies to obtain these values considering  $\bar{p}=1$ and $\bar{p}=\infty$, respectively.
  To do so, we will obtain bounds for the $w$-variables that will be included in the model to strengthen the formulation.

\subsubsection{Initial big M parameters} \label{ss:q1pinf}

Note that  by Proposition \ref{prop1}, using $\bar{q}=1$ and  $\bar{p}=\infty$, we can consider $M_i= \mbox{dist}_i^{1} \cdot\mbox{UB}_{\tiny{\mbox{RL-}\ell_1}}$ as the initial value of $M_i$, for $i\in N$, 
where $\mbox{dist}_i^1 =\underset{j \in N}{\max } \; \{ \norm{x_i-x_j}_1: y_i=y_j \}$ and $\mbox{UB}_{\tiny{\mbox{RL-}\ell_1}}$ is an upper bound of (RL-$\ell_1$-M).  Since $\norm{w}_{\infty}\leq \norm{w}_1$ and $\mbox{UB}_{\tiny{\mbox{RL-}\ell_1}}$ is an upper bound of $\norm{w}_1,$ we have that $\norm{w}_{\infty}\leq \norm{w}_{1} \leq$ $\mbox{UB}_{\tiny{\mbox{RL-}\ell_1}}$. Furthermore, an upper bound of (RL-$\ell_1$-M) can be easily obtained from a feasible solution $({\tilde{w}^+},{\tilde{w}^-}, \tilde{\bb}, \tilde{\xi}, \tilde{z}),$ built from the optimal solution of  (SVM-$\ell_1$), ($w^{\tiny{{\mbox{SVM}}}},\bb^{\tiny{{\mbox{SVM}}}}, \xi^{\tiny{{\mbox{SVM}}}}$), establishing that $\bb=\bb^{\tiny{{\mbox{SVM}}}}$,
\begin{eqnarray*}
 & \tilde{w}^+_k = \begin{cases}
w^{\tiny{{\mbox{SVM}}}}_k, &\text{if } w^{\tiny{{\mbox{SVM}}}}_k\geq 0,\\
0, &\text{otherwise,}
\end{cases}\text{for }k \in D,\;\;\;
&  \tilde{w}^-_k = \begin{cases}
-w^{\tiny{{\mbox{SVM}}}}_k,&\text{if } w^{\tiny{{\mbox{SVM}}}}_k\leq 0,\\
0, &\text{otherwise,}
\end{cases}\text{for }k \in D,\;\;\;
\\
& \tilde{\xi}_i= \begin{cases}
\xi^{\tiny{{\mbox{SVM}}}}_i, &\text{if }\xi^{\tiny{{\mbox{SVM}}}}_i \leq 2,\\
0, &\text{otherwise,}
\end{cases}\text{for }i \in N,\;\;\;
&  \tilde{z}_i = \begin{cases}
0, &\text{if } \xi^{\tiny{{\mbox{SVM}}}}_i \leq 2,\\
1, &\text{otherwise,}
\end{cases} \text{for }i \in N.
\end{eqnarray*}
This upper bound could be improved using the information given by $\tilde{z}$ values to obtain the following model:
\begin{alignat}{4}
(\overline{\mbox{SVM-$\ell_1$}})_{\tilde{z}} & \quad &  & \mbox{min}&  \quad & 	\displaystyle \sum_{k=1}^{d}(w_k^+ + w_k^-) + C  \left( \sum_{i\in N:\tilde{z}_i=0}\xi_i \right), & & \nonumber \\
& & &\mbox{s.t.} & \quad & \eqref{ec:w}, \nonumber \\
& & & & &  y_i \left(\sum_{k=1}^{d}(w_k^+-w_k^-) x_{ik}+\bb \right) \geq 1- \xi_i,     & \quad &i \in N :\tilde{z}_i=0, \nonumber \\
& & & & & 0 \leq   \xi_i \leq 2, & \quad  & i \in N:\tilde{z}_i=0 . \nonumber  
\end{alignat}

The solution of this linear problem $({\bar{w}^+},{\bar{w}^-}, \bar{\bb}, \bar{\xi})$ together with $\tilde{z}$ values constitute a feasible solution for (RL-$\ell_1$-M) that provides a better upper bound, $\mbox{UB}_{\tiny{\mbox{RL-}\ell_1}}$. With this new upper bound, the value of $M_i$ can be updated.

In Variant 1 of Algorithm \ref{st:1a}, we give a pseudocode  with a strategy to obtain a valid value of $M_i$, for $i \in N$. The main purpose of this algorithm is to solve a set of linear problems in order to compute a tightened upper bound of $\norm{w}_\infty$ (Steps 5-7).  Consequentially, we obtain bounds for the $w$-variables and add them to the problem, thus improving the formulation (Step 8). 

On the other hand, if we use $\bar{q}=\infty$ and $\bar{p}=1$, the initial value of $M_i$ for $i\in N$  would be $\left(\underset{j \in N}{\max } \; \{ \norm{x_i-x_j}_{\infty} : y_i=y_j \} \right) \mbox{UB}_{\tiny{\mbox{RL-}\ell_1}},$ because $\mbox{UB}_{\tiny{\mbox{RL-}\ell_1}}$ is an upper bound of $\norm{w}_1$. In order to improve this bound, we propose a strategy which is summarized in Variant 2 of Algorithm \ref{st:1a}, the first step of which is to compute the initial values for the big M parameters. Next,  
we compute a tightened bound of $\norm{w}_1$, thus  solving a linear programming problem (Step 10 of Algorithm \ref{st:1a}). 

\begin{algorithm}[htbp]\label{st:1a}
	\KwData{Training sample composed by a set of $n$ individuals with $d$ features.}
	\KwResult{Update values of $M_i$ and bounds for $w_k^++w_k^-.$ }

	Solve the problem (SVM-$\ell_1$). From its optimal solution, build a feasible  solution of (RL-$\ell_1$-M), $({\tilde{w}^+},{\tilde{w}^-}, \tilde{\bb}, \tilde{\xi}, \tilde{z})$.  Solve ($\overline{\mbox{SVM-}\ell_1}$)$_{\tilde{z}}$ and build an improved feasible solution. Update the upper bound 	$\mbox{UB}_{\tiny{\mbox{RL-}\ell_1}}$.
		
	\For{$i\in N$}{ 
	$\mbox{dist}_i^\infty =\underset{j \in N}{\max } \;\{ \norm{x_i-x_j}_ \infty : y_i=y_j \}, \;\;
	M_i= \mbox{dist}_i^{\infty} \cdot \mbox{UB}_{\tiny{\mbox{RL-}\ell_1}}.$  }
	
	\Case{\text{Variant 1}}{ 
			\For{$k_0 \in D$}{ Solve the following linear programming problem: 
			\begin{alignat}{3}
			& \mbox{max} &  \quad & 
			\displaystyle w_{k_0}^+ + w_{k_0}^-, &&\nonumber \\
			&\mbox{s.t.} &  \quad & \eqref{rl-l2-ec2},  \eqref{rest:cond1}, \eqref{ec:w},\eqref{ec:linealizado},  & \quad & \nonumber \\
			& & & \displaystyle\sum_{k=1}^{d}(w_k^+ + w_k^-) + C  \left( \sum_{i=1}^n\xi_i + 2 \sum_{i=1}^n z_i \right)\leq \mbox{UB}_{\tiny{\mbox{RL-}\ell_1}}, &  & \label{ec:UB} \\
			& & & 0 \leq z_i\leq 1, &   \quad  & i \in N.   \label{ec:rel}
			\end{alignat}
			
			Let $\mbox{UB}_{w_{k_0}}$ be the optimal objective value of the above problem. 
		}
		
		Update $M_i= \min\left\{\mbox{dist}_i^{1} \cdot\underset{k\in D}{ \max } \; \left\{\mbox{UB}_{w_{k}}\right\},\mbox{dist}_i^{\infty} \cdot \mbox{UB}_{\tiny{\mbox{RL-}\ell_1}}\right\}$ and add the obtained bounds to the problem (RL-$\ell_1$-M) including the following set of constraints: 
		\begin{equation}
		w_k^+ + w_k^- \leq \mbox{UB}_{w_{k}}, \quad k \in D. \label{ec:boundwk}
		\end{equation}
		}	
		\Case{\text{Variant 2}}{
				Solve the following linear programming problem: 
			\begin{alignat}{3}
			& \max &  \quad & 
			\displaystyle \sum_{k=1}^{d}w_k^+ + w_k^-, &&\nonumber \\
			&\mbox{s.t.} &  \quad &  \eqref{rl-l2-ec2}, \eqref{rest:cond1},\eqref{ec:w}- \eqref{ec:rel}.  & \quad & \nonumber 
			\end{alignat}
			Let $\mbox{UB}_{w}$ be the optimal objective value of the above problem. 
			
			Update $M_i= \mbox{dist}_i^\infty \cdot \mbox{UB}_{w}$ and add the obtained bounds to the problem (RL-$\ell_1$-M) including the following set of constraints: 
			\begin{equation}
			w_k^+ + w_k^- \leq \mbox{UB}_{w}, \quad k \in D. \label{ec:boundw}
			\end{equation}
		}
			\caption{Variant 1 and 2. Computation of tightened bounds of $w$-variables.}
\end{algorithm}

Observe that in both variants of Algorithm \ref{st:1a}, the initial values of the big M parameters are computed as  $M_i= \mbox{dist}_i^\infty  \cdot \mbox{UB}_{\tiny{\mbox{RL-}\ell_1}}$, where  $\mbox{dist}_i^\infty= \underset{j \in N}{\max } \; \{ \norm{x_i-x_j}_{\infty} : y_i=y_j \},$ because $\mbox{dist}_i^\infty\leq \mbox{dist}_i^1$ as well as $\mbox{UB}_{\tiny{\mbox{RL-}\ell_1}}$ is an upper bound of $\norm{w}_{\infty}$ and $\norm{w}_1$.

A detailed comparison between the performance of both strategies is carried out in Section \ref{sec:computationalExp}. In general, the resolution time of this algorithm will be lower when using  Variant 2 because it solves fewer problems in each iteration.

\subsubsection{A Lagrangian relaxation based procedure}
\label{subsec:Lagrangian-relaxation}

Similarly to the previous strategy, this subsection aims to 
tighten the bounds for the $w$-variables. Up to this point, the bounds on the $w$-variables 
 were computed by solving linear problems and expressed as valid inequalities, i.e., constraints \eqref{ec:boundwk} and \eqref{ec:boundw}, using Variants 1 and 2 respectively. In the following, we will refer to set of constraints \eqref{ec:boundwk}, but the same theoretical results can be obtained using Variant 2, i.e., including \eqref{ec:boundw} instead of \eqref{ec:boundwk}. In this subsection, we propose a procedure that is based on the Lagragian relaxation of  (RL-$\ell_1$-M). Before going into detail about this strategy, note that the linear relaxation of (RL-$\ell_1$-M) is
\begin{alignat}{4}
\mbox{(LP-RL-$\ell_1$)} & \quad &  & \mbox{min}&  \quad & 	\displaystyle \sum_{k=1}^{d}(w_k^+ + w_k^-) + C  \left( \sum_{i=1}^n\xi_i + 2 \sum_{i=1}^n z_i \right), & & \nonumber \\
& & &\mbox{s.t.} &  \quad & \eqref{rl-l2-ec2}, \eqref{rest:cond1},\eqref{ec:w},\eqref{ec:linealizado},\eqref{ec:rel}, \eqref{ec:boundwk}.  & \quad & \nonumber 
\end{alignat}

We will build two new models based on the previous one, seeing as this will be necessary for the next technical result that will be used in Theorem \ref{coro:lg} to provide bounds on the $w$-variables using Lagrangian relaxations. Given a value of $\tilde{w}^+_{k_0}>0,$ for $k_0 \in D,$ we obtain the following equivalent model to (LP-RL-$\ell_1$) by making the changes of variables: $\bar{w}_{k_0}^+= w_{k_0}^+-\tilde{w}^+_{k_0},$ $\bar{w}_{k}^+= w_{k}^+,$ for $k\in D\setminus\{k_0\},$ and $\bar{w}_{k}^-= w_{k}^-,$ for $k\in D:$ 
\begin{alignat}{4}
(\overline{\mbox{LP-RL-$\ell_1$}})_{k_0}^+  & \; \;& & \min&  \; & 	\displaystyle \sum_{k=1}^{d}(\bar{w}_k^+ + \bar{w}_k^-) + C  \left( \sum_{i=1}^n\bar{\xi}_i + 2 \sum_{i=1}^n \bar{z}_i \right)+\tilde{w}^+_{k_0} , & & \nonumber \\
&& &\mbox{s.t.} &   & y_i\left( \sum_{k=1}^{d}(\bar{w}_k^+ - \bar{w}_k^-)  x_{ik}+\bar{\bb}\right) \geq 1- \bar{\xi}_i -y_i\tilde{w}^+_{k_0} x_{ik_0} -M_i\bar{z}_i,\label{ec:linear-wk}  & \quad  & i \in N,  \\
&&& & & \bar{\xi}_i \leq 2(1-\bar{z}_i),  & \quad  & i\in N, \label{ec:cond1prime}\\
&&& & & \bar{w}_k^+ + \bar{w}_k^- \leq \mbox{UB}_{w_{k}}, \quad& &k \in D,\\
 &&& & &  0 \leq   \bar{\xi}_i \leq 2,  & \quad  & i\in N, \label{ec:eprime}\\
 &&& & & 0 \leq \bar{z}_i\leq 1, &   \quad  & i\in N, \label{ec:zprime}\\
 &&& & &  \bar{w}_k^+\geq 0, & \quad &  k\in D\setminus\{k_0\},  \label{ec:wprime} \\
  &&& & &  \bar{w}_{k_0}^+\geq -\tilde{w}^+_{k_0}, \\
  &&& & &\bar{w}_{k}^- \geq 0,  &  & k \in D. \label{ec:wmenos} 
\end{alignat}

Similarly, given a value of $\tilde{w}^-_{k_0}>0,$ for $k_0 \in D,$ the following equivalent model to (LP-RL-$\ell_1$) is obtained by changing of variables $\bar{w}_{k_0}^-= w_{k_0}^--\tilde{w}^-_{k_0},$ $\bar{w}_{k}^+= w_{k}^+,$ for $k\in D,$ and $\bar{w}_{k}^-= w_{k}^-,$ for $k\in D\setminus\{k_0\}$: 
\begin{alignat}{4}
(\overline{\mbox{LP-RL-$\ell_1$}})_{k_0}^- & \;\; &  & \mbox{min}&  \;\; & 	\displaystyle \sum_{k=1}^{d}(\bar{w}_k^+ + \bar{w}_k^-) + C  \left( \sum_{i=1}^n\bar{\xi}_i + 2 \sum_{i=1}^n \bar{z}_i \right)+\tilde{w}^-_{k_0}, & & \nonumber \\
& & &\mbox{s.t.} &  \;\; &  \eqref{ec:cond1prime}-\eqref{ec:zprime},   \nonumber\\
&&& & &  y_i\left( \sum_{k=1}^{d}(\bar{w}_k^+ - \bar{w}_k^-) x_{ik}+\bar{\bb}\right) \geq 1- \bar{\xi}_i +y_i\tilde{w}^-_{k_0} x_{ik_0} -M_i\bar{z}_i,\label{ec:linear-wk-}  & \quad  & i \in N, \\
  &&& & &\bar{w}_{k}^+ \geq 0, &  & k\in D, \label{ec:wmas} \\
   &&& & &\bar{w}_{k}^- \geq 0, &  &  k\in D\setminus\{k_0\}, \\
&&& & &  \bar{w}_{k_0}^-\geq -\tilde{w}^-_{k_0}. \end{alignat}

\begin{lema}\label{tm:lagrange2}
 The following statements hold:
\begin{itemize}
\item[i)] Let ($\bar{w}^{+^\ast}, \bar{w}^{-^\ast}, \bar{\bb}^{\ast}, \bar{\xi}^{\ast}, \bar{z}^{\ast}$) be an optimal solution of $(\overline{\mbox{LP-RL-$\ell_1$}})_{k_0}^+$ with $\bar{w}^{+^\ast}_{k_0}=0$ for $k_0 \in D,$ with $Z_{k_0}^{+}$ being its objective value, and $\bar{\alpha}^+$ a vector of optimal values for the dual variables associated with constraints \eqref{ec:linear-wk}.  If ($\bar{w}^{+^\prime}, \bar{w}^{-^\prime}, \bar{\bb}^{\prime}, \bar{\xi}^{\prime}, \bar{z}^{\prime}$)  is an optimal solution of  $(\overline{\mbox{LP-RL-$\ell_1$}})_{k_0}^+$ restricting $\bar{w}^+_{k_0}= \hat{w}^{+}_{k_0},$ where $\hat{w}^{+}_{k_0} > -\tilde{w}^+_{k_0},$ and $Z_{\hat{k}_0}^{+}$ its objective value, then:  
		\begin{equation}
		Z_{k_0}^{+} +\hat{w}^{+}_{k_0}\left(1- \sum_{i=1}^{n}\bar{\alpha}_i^+ y_i x_{ik_0}\right)\leq Z_{\hat{k}_0}^{+}. 
		\end{equation}
\item[ii)] Let ($\bar{w}^{+^\ast}, \bar{w}^{-^\ast}, \bar{\bb}^{\ast}, \bar{\xi}^{\ast}, \bar{z}^{\ast}$) be an optimal solution of $(\overline{\mbox{LP-RL-$\ell_1$}})_{k_0}^-$  with $\bar{w}^{-^\ast}_{k_0}=0$, for $k_0 \in D,$ with $Z_{k_0}^{-}$ being its objective value, and $\bar{\alpha}^-$ a vector of optimal values for the dual variables associated with constraints \eqref{ec:linear-wk-}.  If ($\bar{w}^{+^\prime}, \bar{w}^{-^\prime}, \bar{\bb}^{\prime}, \bar{\xi}^{\prime}, \bar{z}^{\prime}$) is an optimal solution of  $(\overline{\mbox{LP-RL-$\ell_1$}})_{k_0}^-$ restricting $\bar{w}^{-}_{k_0}= \hat{w}^{-}_{k_0},$ where $\hat{w}^{-}_{k_0}>-\tilde{w}^-_{k_0},$ and $Z_{\hat{k}_0}^-$ its objective value, then:
\begin{equation}
Z_{k_0}^{-} +\hat{w}^{-}_{k_0}\left(1+ \sum_{i=1}^{n}\bar{\alpha}_i^- y_i x_{ik_0}\right)\leq  Z_{\hat{k}_0}^-. 
\end{equation}

\end{itemize}
\end{lema}

\dem
 i) Let $\bar{\alpha}^+$ be the vector of optimal values for the dual variables associated with family of constraints $\eqref{ec:linear-wk}$ for $(\overline{\mbox{LP-RL-$\ell_1$}})_{k_0}^+$. By the complementary slackness conditions, it holds that:
\begin{align}
Z_{k_0}^{+}=& \sum_{k=1}^{d}\left(\bar{w}_k^{+^\ast} + \bar{w}_k^{-^\ast}\right) + C \left(\sum_{i=1}^{n}\bar{\xi}_i^\ast + 2\sum_{i=1}^{n} \bar{z}_i^\ast\right)+\tilde{w}^+_{k_0} \nonumber\\& + \sum_{i=1}^{n} \bar{\alpha}_i^+\left( 1- \bar{\xi}_i^\ast  -y_i\tilde{w}^+_{k_0} x_{ik_0} -M_i\bar{z}_i^\ast- y_i \sum_{k=1}^{d}\left(\bar{w}_k^{+^\ast} - \bar{w}_k^{-^\ast}\right)x_{ik}- y_i\bar{\bb}^\ast\right). \nonumber
\end{align}

Since $\bar{w}_{k_0}^{+^\ast}=0,$ this variable can be removed from the summations, giving:
\begin{align}
Z_{k_0}^{+}=& \sum_{k=1,k\neq k_0}^{d}\bar{w}_k^{+^\ast} +\sum_{k=1}^{d}\bar{w}_k^{-^\ast} + C \left(\sum_{i=1}^{n}\bar{\xi}_i^\ast + 2\sum_{i=1}^{n} \bar{z}_i^\ast\right)+\tilde{w}^+_{k_0} \nonumber\\& + \sum_{i=1}^{n} \bar{\alpha}_i^+\left( 1- \bar{\xi}_i^\ast -y_i\tilde{w}^+_{k_0} x_{ik_0} -M_i\bar{z}_i^\ast- y_i \left(\sum_{k=1,k\neq k_0}^{d}\bar{w}_k^{+^\ast}x_{ik} - \sum_{k=1}^{d}\bar{w}_k^{-^\ast}x_{ik} \right)- y_i\bar{\bb}^\ast\right).\label{cond:lang2}
\end{align}
 
Alternatively, the $(\overline{\mbox{LP-RL-$\ell_1$}})_{k_0}^+$ model with additional constraint $\bar{w}^+_{k_0}= \hat{w}^{+}_{k_0},$ and in which the family of constraints \eqref{ec:linear-wk} has been dualized, with  $\alpha_i\geq 0$ for any $i\in N$, is the following: 
\begin{alignat}{3}
& \mbox{min}&  \quad & 	\displaystyle \sum_{k=1}^{d}\left(\bar{w}_k^{+} + \bar{w}_k^{-}\right)  + C  \left( \sum_{i=1}^n\bar{\xi}_i + 2 \sum_{i=1}^n \bar{z}_i \right) +\tilde{w}^+_{k_0}  \nonumber&&\\
&&& + \sum_{i=1}^{n} \alpha_i\left( 1- \bar{\xi}_i -y_i\tilde{w}^+_{k_0} x_{ik_0}-M_i\bar{z}_i- y_i \sum_{k=1}^{d}\left(\bar{w}_k^{+} - \bar{w}_k^{-}\right)x_{ik}- y_i\bar{\bb}_i\right),   \nonumber&& \\
&\mbox{s.t.} &  \quad & \eqref{ec:cond1prime}-\eqref{ec:wprime},\eqref{ec:wmenos},  && \nonumber \\
& & & \bar{w}^+_{k_0}= \hat{w}^{+}_{k_0}.\nonumber   && 
\end{alignat}
Thus, by extracting the coefficients of $\hat{w}_{k_0}^{+}$ from the summations, this problem can be rewritten as follows: 
\begin{alignat}{4}
(\overline{\mbox{Lg-RL-$\ell_1$}})_{k_0}^+& &\mbox{min} &   	\displaystyle  \sum_{k=1,k\neq k_0}^{d}\bar{w}_k^{+} +\sum_{k=1}^{d}\bar{w}_k^{-} + C  \left( \sum_{i=1}^n\bar{\xi}_i + 2 \sum_{i=1}^n \bar{z}_i \right) +\tilde{w}^+_{k_0} + \hat{w}_{k_0}^{+}\left(1- \sum_{i=1}^{n} \alpha_iy_ix_{ik_0} \right)  + &&\nonumber\\ &&& \sum_{i=1}^{n} \alpha_i\left( 1- \bar{\xi}_i -y_i\tilde{w}^+_{k_0}x_{ik_0} -M_i\bar{z}_i- y_i \left(\sum_{k=1,k\neq k_0}^{d}\bar{w}_k^{+}x_{ik} - \sum_{k=1}^{d}\bar{w}_k^{-}x_{ik} \right)- y_i\bar{\bb}_i\right), & & \nonumber \\
& & \mbox{s.t.} \quad&  \eqref{ec:cond1prime}-
\eqref{ec:wprime},\eqref{ec:wmenos}. \nonumber 
\end{alignat}

Observe that  ($\bar{w}^{+^\ast}, \bar{w}^{-^\ast}, \bar{\bb}^{\ast}, \bar{\xi}^{\ast}, \bar{z}^{\ast}$), an optimal solution of $(\overline{\mbox{LP-RL-$\ell_1$}})_{k_0}^+$, is feasible for the problem above, since all constraints of problem $(\overline{\mbox{Lg-RL-$\ell_1$}})_{k_0}^+$ are included in the former.  
Moreover, any feasible solution of problem $(\overline{\mbox{Lg-RL-$\ell_1$}})_{k_0}^+,$ taking $\bar{w}^+_{k_0}=0,$ is feasible for  $(\overline{\mbox{LP-RL-$\ell_1$}})_{k_0}^+$ where family of constraints \eqref{ec:linear-wk}  has been dualized. Therefore, for $\alpha= \bar{\alpha}^+$, using \eqref{cond:lang2}, the optimal objective value of model $(\overline{\mbox{Lg-RL-$\ell_1$}})_{k_0}^+$ is $Z_{k_0}^{+}+ \hat{w}^{+}_{k_0}\left(1- \displaystyle\sum_{i=1}^{n}\bar{\alpha}_i^+y_ix_{ik_0}\right).$ Since $(\overline{\mbox{Lg-RL-$\ell_1$}})_{k_0}^+$ is a Lagrangian relaxation of $(\overline{\mbox{LP-RL-$\ell_1$}})_{k_0}^+$ with the additional constraint $\bar{w}^+_{k_0}= \hat{w}^{+}_{k_0}$, then the optimal objective value of  $(\overline{\mbox{Lg-RL-$\ell_1$}})_{k_0}^+$ is a lower bound 
of $Z_{\hat{k}_0}^{+}$. Result ii) is proven by following the same argument as i).

\fin

By applying the lemma above, we obtain the following theorem which provides the bounds of the $w$-variables.

\begin{teor}\label{coro:lg}
 
Under the hypothesis of Lemma \ref{tm:lagrange2}
we obtain the following bounds of the $w$-variables for the (RL-$\ell_1$-M) problem: 
	\begin{itemize}
		\item $w_{k_0}^+ \leq  \mbox{UB}_{w_{k}^+}:= \dfrac{\mbox{UB}_{\tiny{\mbox{RL-}\ell_1}}-Z_{k_0}^{+}}{1- \sum_{i=1}^{n} \bar{\alpha}_i^+y_ix_{ik_0} }+\tilde{w}^+_{k_0},$
		\item  $w_{k_0}^- \leq  \mbox{UB}_{w_{k}^-}:= \dfrac{\mbox{UB}_{\tiny{\mbox{RL-}\ell_1}}-Z_{k_0}^{-}}{1+ \sum_{i=1}^{n} \bar{\alpha}_i^-y_ix_{ik_0} }+\tilde{w}^-_{k_0},$
	\end{itemize}

where $\mbox{UB}_{\tiny{\mbox{RL-}\ell_1}}$ is an upper bound of (RL-$\ell_1$-M).
\end{teor}
\dem
An equivalent model $(\overline{\mbox{LP-RL-$\ell_1$}})_{k_0}^+$  $\left((\overline{\mbox{LP-RL-$\ell_1$}})_{k_0}^-)\right)$ can be built from an optimal solution of (LP-RL-$\ell_1$), i.e., $({w^+}^*,{w^-}^*, \bb^*, \xi^*, z^*)$  that satisfies that ${w^+_{k_0}}^*=\tilde{w}^+_{k_0}$ $\left({w^-_{k_0}}^*=\tilde{w}^-_{k_0}\right).$ In this situation, both models have the same optimal objective value. Since model (LP-RL-$\ell_1$) is the linear relaxation of model (RL-$\ell_1$-M), an upper bound of the latter would be an upper bound of former. Thusly, if the objective value of an optimal solution of  $(\overline{\mbox{LP-RL-$\ell_1$}})_{k_0}^+$ restricting $\bar{w}^{+}_{k_0}= \hat{w}^{+}_{k_0}$ $\left((\overline{\mbox{LP-RL-$\ell_1$}})_{k_0}^-\right.$ restricting $\left.\bar{w}^{-}_{k_0}= \hat{w}^{-}_{k_0}\right)$ is bigger than $\mbox{UB}_{\tiny{\mbox{RL-}\ell_1}},$ the value $\bar{w}^{+}_{k_0}= \hat{w}^{+}_{k_0}$ $\left(\bar{w}^{-}_{k_0}= \hat{w}^{-}_{k_0}\right)$ can be discarded as an optimal solution of (RL-$\ell_1$-M). This is because any solution with this value will provide a solution with an objective value that is worse than $\mbox{UB}_{\tiny{\mbox{RL-}\ell_1}}$. Therefore, we can restrict ourselves to the values of $\hat{w}^{+}_{k_0}$ $\left(\hat{w}^{-}_{k_0}\right)$ in such a way that $Z_{\hat{k}_0}^+ \leq \mbox{UB}_{\tiny{\mbox{RL-}\ell_1}}$ $\left(Z_{\hat{k}_0}^- \leq \mbox{UB}_{\tiny{\mbox{RL-}\ell_1}}\right)$. Thus, according to Lemma \ref{tm:lagrange2}, $\bar{w}^{+}_{k_0}$ $\left(\bar{w}^{-}_{k_0}\right)$ satisfies: $$\bar{w}^{+}_{k_0} \leq \dfrac{\mbox{UB}_{\tiny{\mbox{RL-}\ell_1}}-Z_{k_0}^{+}}{1- \sum_{i=1}^{n} \bar{\alpha}_i^+y_ix_{ik_0} }, \quad \bar{w}^{-}_{k_0} \leq \dfrac{\mbox{UB}_{\tiny{\mbox{RL-}\ell_1}}-Z_{k_0}^{+}}{1+ \sum_{i=1}^{n} \bar{\alpha}_i^+y_ix_{ik_0}}.$$
We therefore obtain the upper bounds for the $w$-variables of problem (RL-$\ell_1$-M) by undoing the changes of the variables are obtained.\fin
 
The resulting bounds will be included in the formulation of the problem as the following set of constraints.
	\begin{eqnarray}
w_k^+  \leq \mbox{UB}_{w_{k}^+},  \quad k\in D, \label{ec:boundw+}\\
w_k^-\leq \mbox{UB}_{w_{k}^-}, \quad k\in D. \label{ec:boundw-}
\end{eqnarray} 

These bounds for the $w$-variables will be used in the next section to obtain tightened values of big M parameters in the models.
\subsection{Tightening values of big M parameters} 

The previous strategies were based on tightening the bounds for the $w$-variables. In this subsection, we will take advantage of these bounds to obtain tightened bounds of the  big M parameters by solving linear programming problems.

 First, we present a result which, for a given $i\in N,$ provides a valid value of the big M parameter when solving a linear problem. Note that in this subsection we will refer to set of constraints \eqref{ec:boundwk}, but the same theoretical results can be obtained using Variant 2 of Algorithm \ref{st:1a}, i.e., including \eqref{ec:boundw} instead.

\begin{prop}\label{prop2}
	If $M_{i}$, for all $i\in N,$ is greater than or equal to the optimal objective value of the following linear problem, provided it is not unbounded:
	\begin{alignat}{4}
&\left(\mbox{UB}_{M_{i}}\right)& \quad&\max&  \quad & 	 1 - \xi_{i}-y_{i} \left(\sum_{k=1}^{d}(w_k^+-w_k^-) x_{ik}+\bb \right),   \nonumber \\
&&&\mbox{s.t.} &  \quad &  \eqref{rl-l2-ec2}, \eqref{rest:cond1},\eqref{ec:w}-
\eqref{ec:boundwk},\eqref{ec:boundw+},\eqref{ec:boundw-}, && \nonumber 
\end{alignat}
then problems (RL-$\ell_1$) and (RL-$\ell_1$-M) are equivalent. 
\end{prop}
\dem
It holds directly from \eqref{ec:linealizado} taking the maximum. 
\fin

In datasets with a large number of individuals, solving a linear problem for each $i \in N$ would be an inefficient computation. 
Thus, we now prove a result that allows us to obtain a valid value for the big M parameter associated with each individual of each class, thereby solving two linear problems.

\begin{prop}\label{prop}
	Problems (RL-$\ell_1$) and (RL-$\ell_1$-M) are equivalent if the following statements hold:
\begin{itemize}
    \item[i)] For all $i\in N,$ when $y_{i}=+1,$ $M_{i}$ is greater than or equal to the optimal objective value of the following linear problem, provided it is not unbounded:
	\begin{alignat}{4}
	&\left(\mbox{UB}_{M_+}\right)& \quad& \max&  \quad & 	 1- \left(\sum_{k=1}^{d}w_k^+ \underline{x}_{+k} -\sum_{k=1}^{d}w_k^- \bar{x}_{+k}+\bb \right),   \nonumber \\
	&&&\mbox{s.t.} &  \quad &  \eqref{rl-l2-ec2},\eqref{rest:cond1},\eqref{ec:w}-
	\eqref{ec:boundwk},\eqref{ec:boundw+},\eqref{ec:boundw-}, && \nonumber 
	\end{alignat}
	where  $\underline{x}_{+k}=\underset{i \in N.}{\min } \; \{ x_{ik}: y_i=+1 \}$ and $\bar{x}_{+k}=\underset{i \in N.}{\max } \; \{ x_{ik}: y_i=+1 \},$ for $k\in D$.
	
	\item[ii)] For all $i\in N,$ when $y_{i}=-1,$ $M_{i}$ is greater than or equal to the optimal objective value of the following linear problem, provided it is not unbounded:
	\begin{alignat}{4}
	&\left(\mbox{UB}_{M_-}\right)& \quad& \max&  \quad & 	 1+ \left(\sum_{k=1}^{d}w_k^+ \bar{x}_{-k} -\sum_{k=1}^{d}w_k^- \underline{x}_{-k}+\bb \right),   \nonumber \\
	&&&\mbox{s.t.} &  \quad &  \eqref{rl-l2-ec2}, \eqref{rest:cond1},\eqref{ec:w}-
	\eqref{ec:boundwk},\eqref{ec:boundw+},\eqref{ec:boundw-}, && \nonumber
	\end{alignat}
	where  $\bar{x}_{-k}=\underset{i \in N.}{\max } \; \{ x_{ik}: y_i=-1 \}$ and $\underline{x}_{-k}=\underset{i \in N}{\min } \; \{ x_{ik}: y_i=-1 \}$, for $k\in D$. 
	\end{itemize}
	\end{prop}

\dem
For each $k\in D$ and $i \in N,$ when $y_{i}=+1,$ the following inequalities are satisfied:  $w_k^+x_{ik} \geq w_k^+ \underline{x}_{+k}$ and $-w_k^-x_{ik} \geq -w_k^- \bar{x}_{+k}.$ Taking the summation in $k,$ we have: 
\begin{equation}1 - \left(\sum_{k=1}^{d}(w_k^+-w_k^-) x_{ik}+\bb \right)\leq 1- \left(\sum_{k=1}^{d}w_k^+ \underline{x}_{+k} -\sum_{k=1}^{d}w_k^- \bar{x}_{+k}+\bb \right),\; \; \mbox{for } i\in N. \label{ec:maxres}\end{equation}

Since $0\leq \xi_{i}\leq2$, for $i \in N$ when $y_{i}=+1,$ we can conclude that:
\begin{equation}1 -\xi_{i}- \left(\sum_{k=1}^{d}(w_k^+-w_k^-) x_{ik}+\bb \right)\leq 1- \left(\sum_{k=1}^{d}w_k^+ \underline{x}_{+k} -\sum_{k=1}^{d}w_k^- \bar{x}_{+k}+\bb \right),\; \; \mbox{for } i\in N. \nonumber\end{equation}

Therefore, the optimal solution of problem $\left(\mbox{UB}_{M_+}\right)$ is an upper bound of problem $\left(\mbox{UB}_{M_{i}}\right)$ for $i \in N$ when $y_{i}=+1$. Similarly, it can be proven that the optimal solution of problem $\left(\mbox{UB}_{M_-}\right)$ is an upper bound of problem $\left(\mbox{UB}_{M_{i}}\right)$ for $i \in N$ when $y_{i}=-1$. Therefore, in accordance with Proposition \ref{prop2} the result holds. \fin

Proposition \ref{prop} provides a strategy to obtain valid values for the big M parameters. Obviously these values will be less tightened than the ones resulting from the application of Proposition \ref{prop2}, but they can be obtained more quickly. In order to find an equilibrium between computational costs and the tightness of the big M values in a dataset that contains a large number of individuals, we propose the strategy which follows. 
First, we cluster individuals of the same class in such a way that individuals in the same cluster have certain similarities. This can be done using the ``$k$-median" or the ``$k$-means" algorithms in each class. Let $C_+,C_-$ be the set of clusters of class 1 and class $-1$ respectively. The number of clusters, i.e., the cardinality of sets $C_+$ and $C_-$ will be a parameter selected by the modeler and it should be adapted to the data. 
Then, for each cluster, $c_+\in C_+$ and $c_- \in C_-,$  we apply the proposition below. The proof thereof has been omitted due to its similarity to Proposition~\ref{prop}. 
 
	\begin{prop}\label{propc}
	Problems (RL-$\ell_1$) and (RL-$\ell_1$-M) are equivalent if the following statements hold: 
	\begin{itemize}
	    \item [i)] For all $i\in c_+$ with $c_+\in C_+,$ $M_{i}$ is greater than or equal to the optimal objective value of the following linear problem, provided it is not unbounded:
		\begin{alignat}{4}
		&\left(\mbox{UB}_{Mc_+}\right)& \quad& \max&  \quad & 	 1- \left(\sum_{k=1}^{d}w_k^+ \underline{x}_{c+k} -\sum_{k=1}^{d}w_k^- \bar{x}_{c+k}+\bb \right),   \nonumber \\
		&&&\mbox{s.t.} &  \quad &  \eqref{rl-l2-ec2},\eqref{rest:cond1},
\eqref{ec:w}-
\eqref{ec:boundwk},\eqref{ec:boundw+},\eqref{ec:boundw-}, && \nonumber 
		\end{alignat}
		where  $\underline{x}_{c+k}=\underset{i \in c_+}{\min } \; \{ x_{ik}\}$ and $\bar{x}_{c+k}=\underset{i \in c_+}{\max } \; \{ x_{ik}\},$ for $k\in D$. 

		\item [ii)]  For all $i\in c_-$ with $c_-\in C_-,$ $M_{i}$ is greater than or equal to the optimal objective value of the following linear problem, provided it is not unbounded:
		\begin{alignat}{4}
		&\left(\mbox{UB}_{Mc_-}\right)& \quad& \max&  \quad & 	 1+ \left(\sum_{k=1}^{d}w_k^+ \bar{x}_{c-k} -\sum_{k=1}^{d}w_k^- \underline{x}_{c-k}+\bb \right),   \nonumber \\
		&&&\mbox{s.t.} &  \quad &  \eqref{rl-l2-ec2}, \eqref{rest:cond1},\eqref{ec:w}-
		\eqref{ec:boundwk},\eqref{ec:boundw+},\eqref{ec:boundw-}, && \nonumber 
		\end{alignat}
		where  $\bar{x}_{c-k}=\underset{i \in c_-}{\max } \; \{ x_{ik} \}$ and $\underline{x}_{c-k}=\underset{i \in c_-}{\min } \; \{ x_{ik} \}$, for $k\in D$. 
		\end{itemize}
		
	\end{prop}

Comparatively tighter values of the big M parameters will be obtained by applying Proposition~\ref{propc} instead of Proposition \ref{prop}, but it solves fewer problems than Proposition~\ref{prop2}. Thus, the best strategy will vary depending on the dataset. We will compare them in Section \ref{sec:computationalExp}.

The results of previous propositions can be summarized in Algorithm \ref{a:unbounded} which is described as a pseudocode. In order to avoid unboundedness in the proposed models, we put forward the strategy that follows. 
One of the two variants of Algorithm \ref{st:1a} is applied in the first step, with the goal of updating the values of the big M parameters and the bounds on the $w$-variables. In the next steps, upper and lower bounds for the $\bb$-variable are
obtained and included in the model. Next, the iterative procedure starts with the solution of problem (LP-RL-$\ell_1$) in Step 4. A new model is built from an optimal solution of this problem and it is solved by applying the required transformations in such a way that its optimal solution verifies the hypothesis of  Theorem \ref{coro:lg}. Hence, Theorem \ref{coro:lg} is applied to obtain new bounds.  After that, the bounds are updated and Steps 2 and 3 are repeated.

The direct procedure is executed next. Said  procedure leads to three variants of the algorithm: Variant I (Steps 5-7) applies Proposition \ref{prop2}, obtaining valid values of the big M parameters and solving a linear problem for each individual of the data set; Variant II (Steps 8-10) applies Proposition \ref{prop}, in which only two linear problems are solved; finally, Variant III (Steps 11-16) applies  Proposition \ref{propc}, in which the number of linear problems solved will be selected by the modeler, depending on the number of clusters (different strategies for subdividing the data into clusters can be followed -- the performance of the ``$k$-median" and the ``$k$-means" algorithms are tested in Section \ref{sec:computationalExp}).  

\begin{algorithm}[htbp]\label{a:unbounded}
	\KwData{Training sample composed by a set of $n$ individuals with $d$ features.}
	\KwResult{Update values of $M_i$, $\mbox{UB}_{w_{k}^+},\mbox{UB}_{w_{k}^-},$ and obtain bounds for $\bb$:  $\mbox{LB}_{\bb}$ and $\mbox{UB}_{\bb}$.} 
Apply Variant 1 or 2 of Algorithm \ref{st:1a}.
	
	Obtain lower $\mbox{LB}_{\bb}$ and upper bounds $\mbox{UB}_{\bb}$ of the $\bb$-variable. For doing that, solve the linear problems that follows. Note that if Variant 2 of Algorithm \ref{st:1a} was used in Step 1, set of constraints \eqref{ec:boundwk} should be replaced by \eqref{ec:boundw}.
		\begin{alignat}{3}
  & \mbox{max/min}  &  \quad & 
	\bb, &&\nonumber \\
 &\mbox{s.t.} &  \quad & \eqref{rl-l2-ec2}, \eqref{rest:cond1}, \eqref{ec:w}-
 \eqref{ec:boundwk}.  & \quad & \nonumber 
	\end{alignat}
	
Include the following constraint in the formulation of the problem: 
\begin{equation}
\mbox{LB}_{\bb} \leq \bb \leq \mbox{UB}_{\bb}. \label{eq:b}
\end{equation}
\Do{an improvement of the bounds is obtained}{

Solve the required problems to apply Theorem \ref{coro:lg} and update $w$-bounds ($\mbox{UB}_{w_{k}^+}$ and $\mbox{UB}_{w_{k}^-}$).  
 Repeat Steps 2 and 3 including constraints \eqref{ec:boundw+} and \eqref{ec:boundw-}.

		\Case{\text{Variant I}}{
					\For{$i\in N$}{ Update $M_{i}$ as the optimal value of the problem $\left(\mbox{UB}_{M_{i}}\right)$+\eqref{eq:b}.
			}
		}
		\Case{\text{Variant II}}{ 
				For $i\in N,$ when $y_i=1,$ update $M_i$ as the optimal value of the problem $\left(\mbox{UB}_{M_+}\right)$+\eqref{eq:b}. \\
			For $i\in N,$ when $y_i=-1,$ update $M_i$  as the optimal value of the problem $\left(\mbox{UB}_{M_-}\right)$+\eqref{eq:b}. 
			
		}
		\Case{\text{Variant III}}{ 
		Cluster the individuals of each class (applying the  $k$-median or the $k$-means algorithm). Let $C_+,C_-$ be the set of clusters of class 1 and class $-1$ respectively.
			 
			\For{$c_+\in C_+$}{Update $M_{i},$ for $i\in c_+,$ as the optimal value of the problem $\left(\mbox{UB}_{Mc_+}\right)$+\eqref{eq:b}.
			}
		\For{$c_-\in C_-$}{Update $M_{i},$ for $i\in c_-,$ as the optimal value of the problem $\left(\mbox{UB}_{Mc_-}\right)$+\eqref{eq:b}.
		}
		
}
}

	\caption{Variant I, II, and III. Computation of tightened values of big M parameters in (RL-$\ell_1$-M).}
\end{algorithm}

\newpage

\section{Strategies for the $\ell_2$-norm case}
\label{sec:l2norm}
In this section, we analyze the ramp loss model whilst considering the  $\ell_2$-norm. We propose strategies to find valid values of the big M parameters, with the objective of enhancing the formulation.

 The problem is formulated as follows: 
\begin{alignat}{4}
\mbox{(RL-$\ell_2$)} & \quad &  & \mbox{min}&  \quad & 	\displaystyle \frac{1}{2} \sum_{k=1}^{d}w_k^2  + C  \left( \sum_{i=1}^n\xi_i + 2 \sum_{i=1}^n z_i \right), & & \nonumber \\
& & &\mbox{s.t.} &  \quad &   \eqref{rl-cons}-
\eqref{rl-l2-ec3}.\nonumber
\end{alignat}

As in the $\ell_1$-norm case, the constraints could be linearized by using a big enough constant $M_i$ in \eqref{rl-cons}. As a result, we obtain the following quadratic model:
\begin{alignat}{4}
\mbox{(RL-$\ell_2$-M)} & \quad &  & \mbox{min}&  \quad & 	\displaystyle \frac{1}{2}\sum_{k=1}^{d}w_k^2  + C  \left( \sum_{i=1}^n\xi_i + 2 \sum_{i=1}^n z_i \right), & & \nonumber \\
& & &\mbox{s.t.} &  \quad &  \eqref{rl-l2-ec2}-\eqref{ec:linealizado_lp}.\nonumber
\end{alignat} 
 
Proposition \ref{prop1}  is still valid for obtaining certain values of the big M parameters for formulation (RL-$\ell_2$-M).
We therefore present a corollary that follows on from Proposition \ref{prop1}. It allows us to compute an initial valid value of the big M parameters. These are based on an upper bound on the $\mbox{UB}_{\tiny{\mbox{RL-}\ell_2}}$ model. An upper bound can be easily built from the optimal solution of problem (SVM-$\ell_2$), denoted by ($\tilde{w},\tilde{b},\tilde{\xi},\tilde{z}$), following a similar procedure to the one described in Section \ref{ss:q1pinf}.

Moreover, this upper bound could be improved by using the information given by $\tilde{z}$ values to obtain the following model:
\begin{alignat}{4}
(\overline{\mbox{SVM-$\ell_2$}})_{\tilde{z}} & \quad &  & \mbox{min}&  \quad & 	\displaystyle \frac{1}{2}\sum_{k=1}^{d}w_k^2 + C  \left( \sum_{i\in N:\tilde{z}_i=0}\xi_i \right), & & \nonumber \\
& & &\mbox{s.t.} &  &  y_i \left(\sum_{k=1}^{d}w_k x_{ik}+\bb \right) \geq 1- \xi_i,     & \quad &i \in N :\tilde{z}_i=0, \nonumber \\
& & & & & 0 \leq   \xi_i \leq 2, & \quad  & i \in N:\tilde{z}_i=0 . \nonumber  
\end{alignat}

The solution of this quadratic problem $(\bar{w}, \bar{\bb}, \bar{\xi})$ together with $\tilde{z}$ values constitute a feasible solution for (RL-$\ell_2$-M) that provides a better upper bound, $\mbox{UB}_{\tiny{\mbox{RL-}\ell_2}}$.

\begin{coro}\label{coro:initialMl2}
Taking the values for the big M parameters such that: 
$$M_i\geq \left(\underset{j \in N}{\max } \; \{ \norm{x_i-x_j}_2 : y_i=y_j \} \right) \sqrt{2\mbox{UB}_{\tiny{\mbox{RL-}\ell_2}}}, \quad i\in N,$$then problems (RL-$\ell_2$) and (RL-$\ell_2$-M) are equivalent. 
\end{coro}

 Similarly to the case of the $\ell_1$-norm, we present a procedure which is based on the Lagragian relaxation to tighten the bounds for the $w$-variables. In order to do that, we define the model where the integrality condition of the $z$-variables has been relaxed, i.e.,
\begin{alignat}{4}
\mbox{(Re-RL-$\ell_2$)} & \quad &  & \mbox{min}&  \quad & 	\displaystyle \frac{1}{2} \sum_{k=1}^{d}w_k^2  + C  \left( \sum_{i=1}^n\xi_i + 2 \sum_{i=1}^n z_i \right), & & \nonumber \\
& & &\mbox{s.t.} &  \quad & \eqref{rl-l2-ec2}, \eqref{rest:cond1},\eqref{ec:linealizado},\eqref{ec:rel}.  & \quad & \nonumber 
\end{alignat}

Hence, we build a  formulation which is based on the previous one, making the following change of the variables: given a value of $\tilde{w}_{k_0}$ for some $k_0\in D,$ $\bar{w}_{k_0}= w_{k_0}-\tilde{w}_{k_0}$ and $\bar{w}_{k}= w_{k},$ for $k\in D\setminus\{k_0\}.$ We therefore obtain the following reformulation of the problem above:
\begin{alignat}{4}
(\overline{\mbox{Re-RL-$\ell_2$}})_{k_0} & \quad &  & \mbox{min}&  \quad & 	\displaystyle \frac{1}{2}\sum_{k=1}^{d}\bar{w}_k^2 + C  \left( \sum_{i=1}^n\bar{\xi}_i + 2 \sum_{i=1}^n \bar{z}_i \right)+\dfrac{1}{2}\tilde{w}_{k_0}^2+\tilde{w}_{k_0} \bar{w}_{k_0} , & & \nonumber \\
& & &\mbox{s.t.} &  \quad &   \eqref{ec:cond1prime}, \eqref{ec:eprime},\eqref{ec:zprime} , \nonumber\\
&&& & &  y_i\left( \sum_{k=1}^{d}\bar{w}_k x_{ik}+\bar{\bb}\right) \geq 1- \bar{\xi}_i -y_i\tilde{w}_{k_0} x_{ik_0} -M_i\bar{z}_i,\label{ec:linearwk}  & \quad  & i \in N.
\end{alignat} 

\begin{teor}\label{tm:lagrangel22}
 Let ($\bar{w}^{\ast}, \bar{\bb}^{\ast}, \bar{\xi}^{\ast}, \bar{z}^{\ast}$) be an optimal solution of $(\overline{\mbox{Re-RL-$\ell_2$}})_{k_0},$ with $\bar{w}_{k_0}^{\ast}=0$, $Z_{k_0}$ being its objective value, and $\bar{\alpha}$ a vector of the optimal values for the dual variables associated with constraints \eqref{ec:linearwk}.  If  ($\bar{w}^{\prime}, \bar{\bb}^\prime, \bar{\xi}^\prime, \bar{z}^\prime$)  is an optimal solution of  $(\overline{\mbox{Re-RL-$\ell_2$}})_{k_0}$ restricting $\bar{w}_{k_0}= \hat{w}_{k_0}$, where $\hat{w}_{k_0}$ is a constant value and $ Z_{\hat{k}_0}$ is its objective value, then:
		\begin{equation}
		Z_{k_0} +\hat{w}_{k_0}\left(\frac{1}{2}\hat{w}_{k_0}+ \tilde{w}_{k_0} - \sum_{i=1}^{n}\bar{\alpha}_i y_i x_{ik_0}\right)\leq   Z_{\hat{k}_0}. 
		\end{equation}
\end{teor}
\dem See Appendix \ref{sec:Appendix}. 

Given the theorem above, we obtain Corollary \ref{coro:l2}, which provides bounds  for the $w$-variables.
\begin{coro}\label{coro:l2}
Under the hypothesis of Theorem \ref{tm:lagrangel22}, we obtain the following bounds of the $w$-variables for the $({\mbox{Re-RL-$\ell_2$}})$ problem:
	\begin{equation} \nonumber
\displaystyle  w_{k_0} \geq \mbox{LB}_{w_{k_0}}:= \sum_{i=1}^{n} \bar{\alpha_i}y_ix_{ik_0} - \sqrt{\displaystyle \left(\tilde{w}_{k_0}-\sum_{i=1}^{n} \bar{\alpha_i}y_ix_{ik_0}  \right)^2 -2 \left( Z_{k_0}-\mbox{UB}_{\tiny{\mbox{RL-}\ell_2}}\right)},
\end{equation}
\begin{equation}  \nonumber
w_{k_0}\leq \mbox{UB}_{w_{k_0}}:= \displaystyle \sum_{i=1}^{n} \bar{\alpha_i}y_ix_{ik_0} + \sqrt{\displaystyle \left(\tilde{w}_{k_0}-\sum_{i=1}^{n} \bar{\alpha_i}y_ix_{ik_0}  \right)^2 -2 \left( Z_{k_0}-\mbox{UB}_{\tiny{\mbox{RL-}\ell_2}}\right)},
\end{equation}
where $\mbox{UB}_{\tiny{\mbox{RL-}\ell_2}}$ is an upper bound of (RL-$\ell_2$-M).
\end{coro}
\dem
See Appendix \ref{sec:Appendix}. 
	
As a consequence of the previous result and having computed the foregoing bounds for the $w$-variables, we include the following set of constraints in the formulation of problem (RL-$\ell_2$-M):
 \begin{equation}
		\mbox{LB}_{w_k} \leq w_k \leq \mbox{UB}_{w_k}, \quad k\in D. \label{ec:wl2}
		\end{equation}
		
In the following, we present a result that allows us to compute valid values for the big M parameters, thus solving a problem for each instance of the dataset.	
	
\begin{prop}\label{propDirectl2}
	Problems (RL-$\ell_2$) and (RL-$\ell_2$-M) are equivalent if, for all $i\in N,$ $M_{i}$ is greater than or equal to the optimal objective value of the following problem, provided it is not unbounded:
	\begin{alignat}{4}
	\left(\mbox{UB}^{\; \ell_2}_{M_{i}}\right) & \quad & & \max&  \quad & 	 1 - \xi_{i}-y_{i} \left(\sum_{k=1}^{d}w_k x_{ik}+\bb \right),   \nonumber \\
	&&&\mbox{s.t.} &  \quad &  \eqref{rl-l2-ec2}, \eqref{ec:linealizado_lp},\eqref{ec:rel}, \eqref{ec:wl2}, && \nonumber \\
	&&&&&\displaystyle \frac{1}{2}\sum_{k=1}^{d}w_k^2 + C  \left( \sum_{\ell=1}^n\xi_{\ell} + 2 \sum_{\ell=1}^n z_{\ell} \right)\leq  \mbox{UB}_{\tiny{\mbox{RL-}\ell_2}}. \label{ec:UBl2}
	\end{alignat}
\end{prop}

\dem
It holds directly from \eqref{ec:linealizado_lp} taking the maximum.
\fin

Although the previous result provides  tightened values of the big M parameters, solving a problem for each individual of the dataset tends to be inefficient.  Hence, given a partition of the individuals, in the following result we present a strategy that allows us to obtain valid values of  the big M parameters for individuals in the same element of that partition.

\begin{prop}\label{prop:Mclusterl2+}
 Problems (RL-$\ell_2$) and (RL-$\ell_2$-M) are equivalent if the following statements hold: 
\begin{itemize}
    \item [i)] For all $i\in c_+,$ where $c_+$ is a subset of the individuals with $y_i=+1$,	$M_{i}$ is greater than or equal to the optimal objective value of the following problem, provided it is not unbounded:
\begin{alignat}{4}
&\left(\mbox{UB}^{\; \ell_2}_{M_{c_+}}\right)& \quad& \max&  \quad & 	 1+ \sum_{k=1}^{d}v_k |\bar{x}_{c_+k}| -\bb \nonumber \\
&&&\mbox{s.t.} &  \quad &  \eqref{rl-l2-ec2},
\eqref{ec:linealizado_lp},\eqref{ec:rel}, \eqref{ec:wl2},\eqref{ec:UBl2}, && \nonumber \\
&&&&& -w_k\leq v_k, && k\in D, \label{ec:vk1}\\
&&&&&    w_k \leq v_k, &&k\in D, \label{ec:vk2} \\
&&&&& v_k \leq \max\{|\mbox{LB}_{w_k}|,|\mbox{UB}_{w_k}|\}, &\quad &k\in D,  \label{ec:vk3}
\end{alignat}
where  $|\bar{x}_{c_+k}|=\underset{i \in c_+.}{\max } \; \{ |x_{ik}|\},$ for $k\in D$. 

 \item [ii)] For all $i\in c_-,$ where $c_-$ is a subset of the individuals with $y_i=-1$,  $M_{i}$ is greater than or equal to the optimal objective value of the following problem, provided it is not unbounded:
	\begin{alignat}{4}
	&\left(\mbox{UB}^{\; \ell_2}_{M_{c_-}}\right)& \quad& \max&  \quad & 	 1+ \sum_{k=1}^{d}v_k |\bar{x}_{c_-k}| +\bb \nonumber \\
	&&&\mbox{s.t.} &  \quad &  \eqref{rl-l2-ec2},\eqref{ec:linealizado_lp}, \eqref{ec:rel}, \eqref{ec:wl2}-\eqref{ec:vk3},
	&& \nonumber 
	\end{alignat}
	where  $|\bar{x}_{c_-k}|=\underset{i \in c_-.}{\max } \; \{ |x_{ik}|\},$ for $k\in D$. 
	\end{itemize} 
\end{prop}
 \dem 
	See Appendix \ref{sec:Appendix}.

All the results of this section can be summarized in a strategy to solve the (RL-$\ell_2$-M) model. Said strategy is described below in Algorithm \ref{algo:l2}. Similarly to Algorithm \ref{a:unbounded}, Algorithm \ref{algo:l2} has three different variants, depending on the procedure chosen to compute the big M parameters. Variant I solves a problem for each individual of the dataset. Variant II solves two problems -- one for the individuals of class $y_i=+1$ and another for the individuals of class $y_i=-1$. In Variant III, the individuals are subdivided into clusters and a problem is solved for each cluster.     

\begin{algorithm}[htbp]\label{algo:l2}
{
	\KwData{Training sample composed by a set of $n$ elements with $d$ features.}
	\KwResult{Update values of $M_i,$  $i\in N$ and obtain bounds for $w_k, \; k\in D$ and $\bb$.}
Solve the problem (SVM-$\ell_2$). From its optimal solution  $(\tilde{w},\tilde{\bb},\tilde{\xi},\tilde{z})$, build a feasible  solution of (RL-$\ell_2$-M). Solve ($\overline{\mbox{SVM-}\ell_2}$)$_{\tilde{z}}$ and  build an improved feasible solution. Update the upper bound 
$\mbox{UB}_{\tiny{\mbox{RL-}\ell_2}}$.
		
		\For{$i\in N$}{ 
			$\mbox{dist}_i^2 =\underset{j \in N}{\max } \;\{ \norm{x_i-x_j}_2 : y_i=y_j \}, \;\;
			M_i= \mbox{dist}_i^2 \cdot \sqrt{2\mbox{UB}_{\tiny{\mbox{RL-}\ell_2}}}.$}
  
		Set  $\mbox{LB}_{w_k}=-\sqrt{2\mbox{UB}_{\tiny{\mbox{RL-}\ell_2}}}$ and $\mbox{UB}_{w_k}=\sqrt{2\mbox{UB}_{\tiny{\mbox{RL-}\ell_2}}}$.
		
	\Do{an improvement of the bounds is obtained}{
			\For{$k_0\in D$}{	Obtain tighter values of $w$-bounds 
				solving the following quadratic problems:
				\vspace{-0.3cm}
				 \begin{alignat}{4}
				& &  & \mbox{max/min}  &  \quad & 
				\displaystyle w_{k_0}, &&\nonumber \\
				& & &\mbox{s.t.} &  \quad & \eqref{rl-l2-ec2},
				\eqref{ec:linealizado_lp},
				\eqref{ec:rel},  \eqref{ec:wl2},\eqref{ec:UBl2}. & \quad & \nonumber  
				\end{alignat}
			From the second iteration, include constraint \eqref{ec:loweupperb} in the previous problem.
}
			 Obtain lower and upper bounds for $\bb$-variable, i.e., $	\mbox{UB}_{\bb}$ and $\mbox{LB}_{\bb}$ solving the following problems:		\vspace{-0.5cm} \begin{alignat}{4}
			  &  &  &  \mbox{max/min}  &  \quad & 
			\bb, &&\nonumber \\
			& & &\mbox{s.t.} &  \quad & \eqref{rl-l2-ec2},\eqref{ec:linealizado_lp},\eqref{ec:rel}, \eqref{ec:wl2}, \eqref{ec:UBl2}.  & \quad & \nonumber 
			\end{alignat}
			
	Include  the following constraint in the formulation of the problem:
			\vspace{-0.3cm}
	\begin{equation}
\mbox{LB}_{\bb} \leq \bb \leq \mbox{UB}_{\bb}. \label{ec:loweupperb}
\vspace{-0.3cm}
	\end{equation} 
	Solve the problem (Re-RL-$\ell_2$) +\eqref{ec:loweupperb}, build $(\overline{\mbox{Re-RL-$\ell_2$}})_{k_0} +\eqref{ec:loweupperb}$  and apply Corollary \ref{coro:l2} to update bounds for the $w$-variables: $\mbox{LB}_{w_k}$ and $\mbox{UB}_{w_k}$.
	
	\Case{\text{Variant I}}{
	\For{$i_0\in N$}{ Solve $\left(\mbox{UB}_{M_{i_0}}^{\;\ell_2}\right)+\eqref{ec:loweupperb}$ and update the big M parameters.
	}
}
\Case{\text{Variant II}}{ Update $M_i$ solving $\left(\mbox{UB}^{\; \ell_2}_{M_{c_+}}\right)+\eqref{ec:loweupperb}$ for $c_+=\{i\in N:y_i=1\}$ and $\left(\mbox{UB}^{\; \ell_2}_{M_{c_-}}\right)+\eqref{ec:loweupperb}$ for $c_-=\{i\in N:y_i=-1\}$.
	
}
\Case{\text{Variant III}}{ 
	Subdivide the individuals of each class in clusters. Let $C_+,C_-$ be the set of clusters of class 1 and class $-1$ respectively.
	
	\For{$c_+\in C_+$}{Update $M_{i},$ for $i\in c_+,$ by solving  $\left(\mbox{UB}_{Mc_+}\right)$+\eqref{ec:loweupperb}.
	}
	\For{$c_-\in C_-$}{Update $M_{i},$ for $i\in c_-,$ by solving 
	$\left(\mbox{UB}_{Mc_-}\right)$+\eqref{ec:loweupperb}.
	}
}		
}
}
\caption{\smaller Variant I, II, and III. Computation of tightened values of big M parameters in (RL-$\ell_2$-M).
}
\end{algorithm}

 Observe that in Algorithm \ref{algo:l2},  unlike Algorithm \ref{a:unbounded}, if we consider clusters containing just one individual of the dataset, the values of the big M parameters obtained in Variant III are not the same as those obtained in Variant I. This is because the $w$-variables are free and the bound is found on the addend $\left|w_kx_{ik}\right|,$ for $i\in N, k\in D.$ However, with the $\ell_1$-norm, the positive and negative values of the $w$-variables are identified with $w_k^+$ and $w_k^-$ respectively. This means that we are able to find the bound for the addend $w_kx_{ik},$ for $i\in N, k\in D.$ In fact, in the tested datasets,  the strategy of Proposition \ref{prop:Mclusterl2+} did not provide good results. Although applying the strategy only took a short amount of time, the improvement on the big M parameters was small and in some cases almost zero. 
 
\section{Computational Experiments}
\label{sec:computationalExp}

In this section, we present the results obtained with our computational experiments.  Specifically, we present a comparison of several solution approaches for solving the models (RL-$\ell_1$), (RL-$\ell_1$-M), (RL-$\ell_2$), and (RL-$\ell_2$-M) using simulated and real-life datasets. 
  
The experiments were conducted on an Intel(R) Xeon(R) W-2135 CPU 3.70 GHz 32 GB RAM computer, using CPLEX 12.7.0. in  Concert Technology C++. As done in \cite{bonami}, due to the presence of big M parameters, the relative MIP tolerance and the integrality tolerance were fixed to zero. The remaining parameters were set to their default values.

\subsection{Data}
The computational experiments were carried out on simulated and real-life datasets. The simulated datasets used in our experiments were Type A and Type B datasets proposed in \cite{31Brooks2011}. These datasets were also used in \cite{carrizosaheuristic} and \cite{bonami}. The latter only included experiments conducted on a challenging
	subset of instances with the following characteristics: $n$=100, $d$=2, Type B. 
We used the Type A and Type B datasets with $n=160,$ i.e., 160 individuals, and $d=2,5,10,$ i.e., $2,5,$ or $10$ features, respectively. Furthermore, the percentage of elements in each class of these instances is $50\%$. We denoted them as: $n$``number of individuals''$d$``number of features''``Type of data''. Thusly, n160d2B signifies that the dataset has 160  individuals, two features, and Type B data.  

The real-life datasets are from the UCI repository \cite{DatasetUCI} and from \cite{ijcnn}.  They are specified in Table \ref{tab:RWDS}, where $n$ is the number of individuals, $d$ is the number of features and the last column states the percentage of elements in each class.  Observe that these datasets have been used to analyze the performance of ramp loss models (see \cite{31Brooks2011, carrizosaheuristic,RLLPSVM}).

\begin{table}[htbp]
	\centering
	\begin{tabular}{llrrr}
		\toprule
		Label & Name in repository  & $n$     & $d$     & Class(\%) \\
		\midrule
		Wpbc   & Breast cancer Wisconsin (prognostic) & 194   & 30    &  76/24\\
		
		SONAR  & Connectionist bench (sonar, mines vs. rocks) & 208   & 60     & 54/46 \\
		
		SPECT   & Cardiac Single Proton Emission Computed Tomography images &  267  &   22  & 79/21 \\
		
		IONO & Ionosphere & 351   & 33    & 64/36 \\
		
		Wdbc   & Breast cancer Wisconsin (diagnostic) & 569   & 30    & 63/37 \\
		WBC   & Breast cancer Wisconsin (original) & 683   & 9    & 65/35 \\
		Ijcnn1 & IJCNN 2001 Neural Network Competition & 35000 & 22    & 91/9 \\
		\bottomrule
	\end{tabular}%
	\caption{Real-life datasets}
	\label{tab:RWDS}%
\end{table}%
In the simulated and real-life datasets, the analysis uses $C \in \{0.01,0.1,1,10,100\}$ as the values for parameter $C,$ as done in \cite{31Brooks2011}, and a time limit of 7200 seconds is imposed, unless stated otherwise. 

\subsection{Exact procedure for (RL-$\pmb{\ell_1}$-M)}

In this section, the exact approaches proposed in Section \ref{sec:Strategies} for solving the \text{(RL-$\ell_1$-M)} formulation are tested and compared with the (RL-$\ell_1$) formulation. We study the performance of different variants using several simulated and real-life datasets. To simplify the notation, we write the name of the algorithm followed by the chosen variant, e.g. Algorithm \ref{a:unbounded}.1.II represents Algorithm \ref{a:unbounded} with Variant 1 in Step 1, i.e., we solve Algorithm \ref{st:1a} using Variant 1 and Variant II (Steps 9 and 10). In each strategy, we also analyze the improvement of the big M parameters. To do this, we compare the values obtained in Step 3 of Algorithm \ref{st:1a}, $M_i^{\text {initial}}$,  with the  tightened values obtained at the end of the applied strategy, $M_i^{\text{final}}$,  i.e., we compute the improvement of M, M's Impr. $=\dfrac{1}{n}\displaystyle\sum_{i=1}^n{\dfrac{M_i^{\text {initial}}-M_i^{\text{final}}}{M_i^{\text {initial}}}}$.  Therefore, this measure depends on $M_i^{\text {initial}}$. As a consequence, we can only compare these values when the strategies are applied to the same data.
  
We first compare two different methods for solving \text{(RL-$\ell_1$-M)} and one for solving (RL-$\ell_1$) on simulated data:  a) Algorithm \ref{a:unbounded}.1.I, b)  Algorithm \ref{a:unbounded}.2.I,  and c) by formulating the model with indicator constraints that generate locally valid implied bound cuts. The latter method was developed by \cite{bonami} and it has been a feature 
of IBM-Cplex since version 12.6.1. We used the method included in 
IBM-Cplex 12.7.0. This solution approach is denoted as Ind. Const. (LIC). 

The results are shown in Table \ref{tab: StIInorm1}, with the first column indicating the name of the dataset and the second column the value of parameter $C$. The next eight columns contain information about strategy performance and the resolution process:  the first and the fifth columns detail the improvement of M during the procedure; the second and the sixth show the strategy time; the third and the seventh indicate the total time, i.e., strategy time plus final time; and the fourth and the eighth report the MIP relative GAP within the time limit. The next groups of columns show information about the behavior of Ind. Const. (LIC): the first column details the total time and the second is the GAP obtained  within the time limit (7200 seconds). When the problem was solved to optimality, the strategy that took the least time is shown in bold, otherwise, the strategy that provided the best GAP within the time limit is highlighted.

\begin{table}[htbp]
	\centering
	\resizebox{\hsize}{!}{
	\begin{tabular}{|c|r|rrrr|rrrr|rr|}
		\hline
		\multicolumn{1}{|l|}{\multirow{2}[1]{*}{Data}} & \multicolumn{1}{r|}{\multirow{2}[1]{*}{$C$}} & \multicolumn{4}{c|}{Algorithm \ref{a:unbounded}.1.I} & \multicolumn{4}{c|}{Algorithm \ref{a:unbounded}.2.I} & \multicolumn{2}{c|}{Ind. Const. (LIC)} \\
		&       & \multicolumn{1}{r}{M's Impr. } & \multicolumn{1}{r}{$t_{st}$} & \multicolumn{1}{r}{$t_{total}$} & \multicolumn{1}{r|}{GAP} & \multicolumn{1}{r}{M's Impr. } & \multicolumn{1}{r}{$t_{st}$} & \multicolumn{1}{r}{$t_{total}$} & \multicolumn{1}{r|}{GAP} & \multicolumn{1}{r}{$t$} & \multicolumn{1}{r|}{GAP} \\
			\hline
\multirow{5}[2]{*}{\begin{turn}{90}n160d2A\end{turn}}  & 100   & \textbf{96.77\%} & \textbf{1.03} & \textbf{2.14} & \textbf{0.00\%} & 96.77\% & 1.19  & 36.19 & 0.00\% & 7207.89 & 52.50\% \\
& 10    & \textbf{96.78\%} & \textbf{1.16} & \textbf{2.12} & \textbf{0.00\%} & 96.78\% & 1.17  & 2.17  & 0.00\% & 7215.59 & 57.82\% \\
& 1     & 96.91\% & 1.47  & 1.67  & 0.00\% & \textbf{97.08\%} & \textbf{1.44} & \textbf{1.64} & \textbf{0.00\%} & 7207.02 & 70.17\% \\
& 0.1   & \textbf{4.19\%} & \textbf{0.63} & \textbf{1.35} & \textbf{0.00\%} & 25.02\% & 1.59  & 5.12  & 0.00\% & 7205.92 & 84.81\% \\
& 0.01  & \textbf{99.95\%} & \textbf{0.11} & \textbf{0.11} & \textbf{0.00\%} & \textbf{99.95\%} & \textbf{0.11} & \textbf{0.11} & \textbf{0.00\%} & 7208.62 & 82.28\%  \\
\hline
\multirow{5}[2]{*}{\begin{turn}{90}n160d5A\end{turn}} & 100   & \textbf{95.09\%} & \textbf{1.20} & \textbf{7202.99} & \textbf{34.95\%} & 95.09\% & 1.11  & 7202.82 & 37.04\% & 7233.20 & 58.56\%  \\
 & 10    & 94.94\% & 1.34  & 7201.46 & 7.91\% & \textbf{94.94\%} & \textbf{1.33} & \textbf{7206.68} & \textbf{5.50\%} & 7206.77 & 64.85\% \\
 & 1     & 94.13\% & 1.58  & 10.02 & 0.00\% & \textbf{94.11\%} & \textbf{1.68} & \textbf{9.11} & \textbf{0.00\%} & 7205.77 & 82.05\% \\
 & 0.1   & 27.69\% & 1.76  & 3.88  & 0.00\% & \textbf{38.04\%} & \textbf{1.91} & \textbf{3.59} & \textbf{0.00\%} & 7214.75 & 73.38\% \\
 & 0.01  & 99.95\% & 0.11  & 0.11  & 0.00\% & \textbf{99.97\%} & \textbf{0.11} & \textbf{0.11} & \textbf{0.00\%} & 7207.12 & 82.05\%  \\
\hline
\multirow{5}[2]{*}{\begin{turn}{90}n160d10A\end{turn}} & 100   & \textbf{93.16\%} & \textbf{1.56} & \textbf{7203.39} & \textbf{52.42\%} & 93.16\% & 1.50  & 7205.10 & 53.09\% & 7204.69 & 69.99\%  \\
 & 10    & 92.89\% & 1.64  & 7216.17 & 44.64\% & \textbf{92.91\%} & \textbf{1.42} & \textbf{7203.62} & \textbf{42.90\%} & 7206.72 & 70.32\% \\
 & 1     & 91.32\% & 2.50  & 7214.95 & 8.78\% & \textbf{91.33\%} & \textbf{2.37} & \textbf{7214.92} & \textbf{7.70\%} & 7219.48 & 77.56\% \\
 & 0.1   & 9.31\% & 1.04  & 74.70 & 0.00\% & \textbf{33.01\%} & \textbf{2.53} & \textbf{5.39} & \textbf{0.00\%} & 7208.65 & 82.05\% \\
 & 0.01  & \textbf{99.80\%} & \textbf{0.13} & \textbf{0.13} & \textbf{0.00\%} & 100.00\% & 0.13  & 0.14  & 0.00\% & 7208.42 & 82.05\%  \\
\hline
\multirow{5}[2]{*}{\begin{turn}{90}n160d2B\end{turn}} & 100   & \textbf{55.55\%} & \textbf{1.16} & \textbf{7203.32} & \textbf{52.15\%} & 55.56\% & 1.09  & 7202.90 & 52.18\% & 7205.64 & 62.71\%  \\
 & 10    & 55.56\% & 1.11  & 7203.02 & 5.90\% & \textbf{55.60\%} & \textbf{1.03} & \textbf{7207.48} & \textbf{4.27\%} & 7213.55 & 65.13\% \\
 & 1     & 55.59\% & 1.33  & 3.28  & 0.00\% & \textbf{56.00\%} & \textbf{1.31} & \textbf{3.23} & \textbf{0.00\%} & 7207.61 & 71.05\% \\
 & 0.1   & 56.57\% & 1.40  & 1.92  & 0.00\% & \textbf{60.85\%} & \textbf{1.39} & \textbf{1.78} & \textbf{0.00\%} & 7205.80 & 81.16\% \\
 & 0.01  & \textbf{96.66\%} & \textbf{0.11} & \textbf{0.11} & \textbf{0.00\%} & 99.93\% & 0.16  & 0.16  & 0.00\% & 7207.97 & 79.71\%  \\
\hline
\multirow{5}[2]{*}{\begin{turn}{90}n160d5B\end{turn}} & 100   & \textbf{72.92\%} & \textbf{1.94} & \textbf{7203.99} & \textbf{58.47\%} & 72.92\% & 1.72  & 7203.87 & 58.98\% & 7206.19 & 69.22\%  \\
 & 10    & \textbf{72.73\%} & \textbf{1.52} & \textbf{7203.48} & \textbf{57.97\%} & 72.73\% & 1.53  & 7204.19 & 58.45\% & 7209.94 & 71.63\% \\
 & 1     & 69.07\% & 1.80  & 7204.08 & 19.55\% & \textbf{69.12\%} & \textbf{1.77} & \textbf{7209.23} & \textbf{16.87\%} & 7207.62 & 72.30\% \\
 & 0.1   & 47.41\% & 1.27  & 3.70  & 0.00\% & \textbf{51.88\%} & \textbf{1.22} & \textbf{3.67} & \textbf{0.00\%} & 7207.81 & 80.57\% \\
 & 0.01  & \textbf{99.95\%} & \textbf{0.11} & \textbf{0.11} & \textbf{0.00\%} & 99.95\% & 0.11  & 0.13  & 0.00\% & 7207.92 & 80.82\%  \\
\hline
\multirow{5}[2]{*}{\begin{turn}{90}n160d10B\end{turn}} & 100   & 60.77\% & 2.10  & 7204.75 & 48.42\% & \textbf{60.77\%} & \textbf{1.88} & \textbf{7203.59} & \textbf{47.42\%} & 7206.53 & 67.00\%  \\
 & 10    & 59.31\% & 2.17  & 7204.57 & 51.76\% & \textbf{59.42\%} & \textbf{2.10} & \textbf{7204.57} & \textbf{48.42\%} & 7208.75 & 58.20\% \\
 & 1     & 79.95\% & 1.78  & 7213.97 & 23.85\% & \textbf{80.12\%} & \textbf{1.78} & \textbf{7213.57} & \textbf{20.96\%} & 7215.72 & 66.13\% \\
 & 0.1   & 42.01\% & 1.19  & 173.63 & 0.00\% & \textbf{45.90\%} & \textbf{1.19} & \textbf{137.97} & \textbf{0.00\%} & 7207.61 & 80.01\% \\
 & 0.01  & 99.98\% & 0.13  & 0.14  & 0.00\% & \textbf{99.79\%} & \textbf{0.13} & \textbf{0.13} & \textbf{0.00\%} & 7208.92 & 80.82\%  \\
\hline

\end{tabular}%
}
	\caption{Performance of exact approaches on simulated data for solving (RL-$\ell_1$-M) and (RL-$\ell_1$).}
	\label{tab: StIInorm1}%

\end{table}%

Note that our strategies improve the  behavior  of Ind. Const. (LIC) in all cases. In fact, in dataset n160d2A, the Ind. Const. (LIC) approach  has a GAP of between 52.50\% and 84.81\% in two hours and, using our strategy, all the individuals are solved to optimality in less than 3 seconds. In the majority of cases, the performance of Algorithm \ref{a:unbounded} using Variant 2 is better than using Variant 1. Also, when Variant 1 works better than Variant 2 there are no significant differences. We can therefore conclude that, when applied to these simulated datasets, the strategy with the best performance is Algorithm \ref{a:unbounded}.2.I. A preliminary test was carried out by solving the problems with the application of Variant II, but the obtained GAPs were worse than with Variant I. Since the time employed in the strategy is less than three seconds in all cases, there is no need to apply Variant III, in which the values of the resulting big M parameters would be less tightened. Note that a $100\%$ of improvement of M does not mean that all the big M parameters are zero, seeing as the improvement of M is the average of the improvement of each $M_i,$ for $i\in N,$ and several of them could be improved by over $100\%$.

Next, we compare the three different methods for solving \text{(RL-$\ell_1$-M)} and (RL-$\ell_1$) respectively on real-life datasets:  a)  Algorithm \ref{a:unbounded}.1.I, b)  Algorithm \ref{a:unbounded}.2.I,  and c) Ind. Const. (LIC). The results are shown in Table \ref{tab: StIInorm1realdata}, which uses a similar notation to Table \ref{tab: StIInorm1}.

\begin{table}[htbp]
	\centering
	\resizebox{0.97\textwidth}{!}{
			\begin{tabular}{|c|r|rrrr|rrrr|rr|}
			\hline
			\multicolumn{1}{|l|}{\multirow{2}[1]{*}{Data}} & \multicolumn{1}{r|}{\multirow{2}[1]{*}{$C$}} & \multicolumn{4}{c|}{Algorithm \ref{a:unbounded}.1.I} & \multicolumn{4}{c|}{Algorithm \ref{a:unbounded}.2.I} & \multicolumn{2}{c|}{Ind. Const. (LIC)} \\
			&       & \multicolumn{1}{r}{M's Impr. } & \multicolumn{1}{r}{$t_{st}$} & \multicolumn{1}{r}{$t_{total}$} & \multicolumn{1}{r|}{GAP} & \multicolumn{1}{r}{M's Impr. } & \multicolumn{1}{r}{$t_{st}$} & \multicolumn{1}{r}{$t_{total}$} & \multicolumn{1}{r|}{GAP} & \multicolumn{1}{r}{$t$} & \multicolumn{1}{r|}{GAP} \\
			\hline
	\multirow{5}[2]{*}{\begin{turn}{90}Wpbc\end{turn}}		  & 100   & 82.88\% & 3.78  & 7206.12 & 63.47\% & \textbf{82.88\%} & \textbf{3.75} & \textbf{7205.97} & \textbf{62.11\%} & 7201.56 & 73.44\% \\
  & 10    & 81.76\% & 4.27  & 7217.15 & 59.72\% & \textbf{81.75\%} & \textbf{4.22} & \textbf{7216.50} & \textbf{58.19\%} & 7203.06 & 72.71\% \\
  & 1     & \textbf{66.11\%} & \textbf{3.54} & \textbf{7216.89} & \textbf{35.94\%} & 66.21\% & 3.41  & 7217.17 & 36.63\% & 7210.61 & 70.74\% \\
  & 0.1   & \textbf{62.62\%} & \textbf{2.32} & \textbf{2.53} & \textbf{0.00\%} & 69.41\% & 2.84  & 3.06  & 0.00\% & 7206.17 & 67.39\% \\
  & 0.01     & \textbf{99.97\%} & \textbf{0.22} & \textbf{0.22} & \textbf{0.00\%} & \textbf{99.94\%} & \textbf{0.22} & \textbf{0.22} & \textbf{0.00\%} & 7205.88 & 67.39\% \\
\hline
	\multirow{5}[2]{*}{\begin{turn}{90}SONAR\end{turn}}	 & 100   & 97.92\% & 11.83 & 12.38 & 0.00\% & 97.92\% & 11.01 & 11.90 & 0.00\% & \textbf{2.20} & \textbf{0.00\%} \\
 & 10    & \textbf{90.84\%} & \textbf{10.25} & \textbf{2058.46} & \textbf{0.00\%} & 90.83\% & 8.76  & 2686.02 & 0.00\% & 7201.62 & 29.98\% \\
 & 1   & 78.13\% & 5.00  & 7217.93 & 46.20\% & \textbf{78.07\%} & \textbf{4.59} & \textbf{7218.52} & \textbf{43.16\%} & 7202.05 & 69.25\% \\
 & 0.1     & 64.74\% & 2.66  & 7215.46 & 9.29\% & \textbf{64.72\%} & \textbf{2.42} & \textbf{7214.13} & \textbf{8.94\%} & 7207.08 & 81.12\% \\
 & 0.01     & 28.91\% & 0.89  & 4.32  & 0.00\% & \textbf{73.80\%} & \textbf{1.53} & \textbf{2.42} & \textbf{0.00\%} & 7206.15 & 85.57\%  \\
\hline
	\multirow{5}[2]{*}{\begin{turn}{90}SPECT\end{turn}}	 & 100   & 63.26\% & 4.82  & 7206.89 & 51.75\% & 63.21\% & 4.83  & 7211.33 & 51.18\% & \textbf{7202.47} & \textbf{34.16\%}  \\
 & 10    & \textbf{64.56\%} & \textbf{6.43} & \textbf{7208.20} & \textbf{27.37\%} & 64.22\% & 4.75  & 7206.52 & 31.43\% & 7203.11 & 36.40\% \\
 & 1     & 60.43\% & 4.92  & 7209.06 & 16.96\% & \textbf{60.44\%} & \textbf{4.95} & \textbf{7207.00} & \textbf{14.51\%} & 7205.19 & 42.86\% \\
 & 0.1   & 53.06\% & 3.87  & 17.63 & 0.00\% & \textbf{53.77\%} & \textbf{4.00} & \textbf{16.31} & \textbf{0.00\%} & 7203.19 & 43.64\% \\
 & 0.01     & 74.25\% & 2.61  & 3.21  & 0.00\% & \textbf{99.55\%} & \textbf{0.28} & \textbf{0.31} & \textbf{0.00\%} & 7204.44 & 32.73\%  \\
\hline
\multirow{5}[2]{*}{\begin{turn}{90}IONO\end{turn}}  & 100   & \textbf{94.95\%} & \textbf{29.86} & \textbf{153.58} & \textbf{0.00\%} & 94.92\% & 29.38 & 159.09 & 0.00\% & 7201.85 & 34.23\%  \\
 & 10    & 90.97\% & 26.77 & 7228.58 & 19.11\% & \textbf{90.97\%} & \textbf{25.61} & \textbf{7227.58} & \textbf{18.93\%} & 7202.60 & 57.00\% \\
 & 1     & \textbf{85.00\%} & \textbf{24.11} & \textbf{7237.08} & \textbf{41.99\%} & 85.06\% & 22.71 & 7236.20 & 42.66\% & 7201.71 & 65.34\% \\
 & 0.1   & \textbf{78.91\%} & \textbf{18.68} & \textbf{7232.26} & \textbf{20.20\%} & 78.95\% & 18.55 & 7232.33 & 22.22\% & 7204.15 & 81.53\% \\
 & 0.01     & \textbf{42.85\%} & \textbf{11.57} & \textbf{16.17} & \textbf{0.00\%} & 77.23\% & 16.76 & 19.15 & 0.00\% & 7205.09 & 88.89\% \\
\hline
\multirow{5}[2]{*}{\begin{turn}{90}Wdbc\end{turn}} & 100   & 99.09\% & 41.83 & 42.50 & 0.00\% & 99.08\% & 38.97 & 39.74 & 0.00\% & \textbf{11.62} & \textbf{0.00\%} \\
  & 10    & 98.18\% & 40.43 & 46.01 & 0.00\% & \textbf{98.12\%} & \textbf{38.37} & \textbf{43.90} & \textbf{0.00\%} & 2432.62 & 0.00\% \\
  & 1     & 101.34\% & 45.64 & 46.00 & 0.00\% & \textbf{101.33\%} & \textbf{40.90} & \textbf{41.18} & \textbf{0.00\%} & 7204.85 & 63.00\% \\
  & 0.1   & 102.21\% & 37.76 & 37.82 & 0.00\% & \textbf{102.22\%} & \textbf{29.47} & \textbf{29.58} & \textbf{0.00\%} & 7202.43 & 86.06\% \\
  & 0.01    & 32.13\% & 23.52 & 27.09 & 0.00\% & \textbf{38.52\%} & \textbf{19.80} & \textbf{23.72} & \textbf{0.00\%} & 7204.90 & 93.87\% \\
\hline
\multirow{5}[2]{*}{\begin{turn}{90}WBC\end{turn}}   & 100   & \textbf{95.88\%} & \textbf{52.02} & \textbf{149.29} & \textbf{0.00\%} & 95.88\% & 50.90 & 163.51 & 0.00\% & 7212.06 & 39.09\% \\
  & 10    & \textbf{95.87\%} & \textbf{43.90} & \textbf{125.14} & \textbf{0.00\%} & 95.78\% & 44.25 & 135.87 & 0.00\% & 7207.31 & 30.46\% \\
  & 1     & \textbf{95.90\%} & \textbf{34.68} & \textbf{180.94} & \textbf{0.00\%} & 95.89\% & 34.13 & 183.21 & 0.00\% & 7207.86 & 39.13\% \\
  & 0.1   & 98.11\% & 32.47 & 33.04 & 0.00\% & \textbf{97.54\%} & \textbf{32.11} & \textbf{32.77} & \textbf{0.00\%} & 7203.27 & 51.23\% \\
  & 0.01     & 97.69\% & 22.79 & 22.85 & 0.00\% & \textbf{97.69\%} & \textbf{21.76} & \textbf{21.82} & \textbf{0.00\%} & 7205.73 & 82.77\%  \\
	
				\hline
			\end{tabular}%

	}
	\caption{Performance of exact approaches on real data for solving (RL-$\ell_1$-M) and (RL-$\ell_1$).}
	\label{tab: StIInorm1realdata}%
	
\end{table}%

The results in Table \ref{tab: StIInorm1realdata}  show  that our strategies improve the behavior of Ind. Const. (LIC)  on real-life datasets in almost all the cases. In fact, when the problem can be solved to optimality in less than two hours with the application of the Ind. Const. (LIC) approach, it can be solved to optimality using our approaches. Conversely, the opposite is not true.  It is worth  highlighting the WBC dataset instances where the Ind. Const. (LIC) approach has a GAP of between 30.46\% and 82.77\% in two hours. However, when using our strategy, all the individuals are solved to optimality in less than three minutes. In the majority of cases, the performance of Algorithm \ref{a:unbounded}.2.I is better than Algorithm \ref{a:unbounded}.1.I.  A preliminary test was carried out by applying Variant II, but the obtained GAPs were worse than when Variant I was used. This is with the exception of the SPECT dataset, the results of which will be analyzed later. Since the strategy is completed in less than one minute in all cases, it is not worth applying Variant III.  

Regarding the time employed in the strategy, note that in the majority of cases Algorithm 2.2.I requires less time than Algorithm 2.1.I. This seems logical because Algorithm 2.2.I solves fewer problems in each iteration. However, there are some cases, e.g. in the IONO dataset in which the value of parameter $C$ is equal to $0.01,$ where the strategy completion time is longer in Algorithm 2.2.I than in Algorithm 2.1.I. This is because the number of iterations is not the same in both procedures. In fact, in the given example (IONO $C=0.01$), the number of iterations of Algorithm 2.1.I, i.e., the number of times that Steps 4-7 are executed, is bigger than the number of iterations of Algorithm 2.2.I.

We will now discuss the particular case of the SPECT dataset.  Table \ref{tab:StIInorm1SPECT} compares the performances of Algorithm \ref{a:unbounded}.1.II and Algorithm \ref{a:unbounded}.2.II with the strategies analyzed previously. The first column states the name of the dataset and the second column the value of parameter $C$.  The table also contains three groups of columns with information about the three different strategies: a) the performance of the best strategy tested in Table \ref{tab: StIInorm1realdata}, b) the behavior of Algorithm \ref{a:unbounded}.1.II, and c) the performance of Algorithm \ref{a:unbounded}.2.II. We can observe that the performance of the previous best strategy from Table \ref{tab: StIInorm1realdata} is worse than using Algorithm \ref{a:unbounded} with Variant II for all the values of $C$ tested.  Although the big M parameters are less tightened using Variant II, the resolution of the problem is faster. If we focus our attention on $C=100, C=10,$ and $C=1$, we note that when using the previous strategies the problems are not solved to optimality in two hours, but when Variant II is used, they are solved to optimality in less than $29, 184,$ and $30$ seconds respectively. In view of these results, we can claim that tighter values of the big M parameters (on average) do not always correspond to a better performance in general, but they do in the majority of the cases. 

\begin{table}[htbp]
	\centering
	\resizebox{0.97\textwidth}{!}{
		\begin{tabular}{|c|r|rrrr|rrrr|rrrr|}
	\hline
	\multicolumn{1}{|l|}{\multirow{2}[2]{*}{Data}} & \multicolumn{1}{r|}{\multirow{2}[2]{*}{$C$}} & \multicolumn{4}{c|}{Best method of Table \ref{tab: StIInorm1realdata} } & \multicolumn{4}{c|}{Algorithm \ref{a:unbounded}.1.II} & \multicolumn{4}{c|}{Algorithm \ref{a:unbounded}.2.II} \\
	&       & \multicolumn{1}{r}{M's Impr. } & \multicolumn{1}{r}{$t_{st}$} & \multicolumn{1}{r}{$t_{total}$} & \multicolumn{1}{r|}{GAP} & \multicolumn{1}{r}{M's Impr. } & \multicolumn{1}{r}{$t_{st}$} & \multicolumn{1}{r}{$t_{total}$} & \multicolumn{1}{r|}{GAP} & \multicolumn{1}{r}{M's Impr. } & \multicolumn{1}{r}{$t_{st}$} & \multicolumn{1}{r}{$t_{total}$} & \multicolumn{1}{r|}{GAP} \\
	\hline

\multirow{5}[2]{*}{\begin{turn}{90}SPECT\end{turn}} & 100   & \multicolumn{1}{l}{-} & \multicolumn{1}{l}{-} & 7202.47 & 34.16\% & \textbf{41.22\%} & \textbf{0.20} & \textbf{28.46} & \textbf{0.00\%} & 41.22\% & 0.16  & 29.98 & 0.00\% \\
 & 10    & 64.56\% & 6.43  & 7208.20 & 27.37\% & 41.42\% & 0.14  & 214.82 & 0.00\% & \textbf{41.38\%} & \textbf{0.16} & \textbf{183.75} & \textbf{0.00\%} \\
 & 1     & 60.44\% & 4.95  & 7207.00 & 14.51\% & \textbf{41.72\%} & \textbf{0.19} & \textbf{29.30} & \textbf{0.00\%} & 41.79\% & 0.16  & 38.12 & 0.00\% \\
 & 0.1   & 53.77\% & 4.00  & 16.31 & 0.00\% & 42.60\% & 0.14  & 3.94  & 0.00\% & \textbf{44.98\%} & \textbf{0.11} & \textbf{3.31} & \textbf{0.00\%} \\
 & 0.01  & 99.55\% & 0.28  & 0.31  & 0.00\% & 75.08\% & 0.13  & 0.70  & 0.00\% & \textbf{97.77\%} & \textbf{0.05} & \textbf{0.08} & \textbf{0.00\%} \\
		\hline

\end{tabular}%
}
	\caption{Performance of Algorithm \ref{a:unbounded} Variant II on SPECT for solving (RL-$\ell_1$-M).}
\label{tab:StIInorm1SPECT}%

\end{table}%

In this paragraph, we analyze the Ijcnn1 dataset which contains a large number of individuals. In particular, we test three random subsets of 2000, 3500, and 5000 individuals respectively. These datasets maintain the same percentage of each class as the original dataset. Additionally, Ijcnn1\_2000 is contained in Ijcnn1\_3500 and this, in turn, is contained in Ijcnn1\_5000. We compare the performance of the  three variants (I, II, and II) of Algorithm \ref{a:unbounded} using Variant 2, because after a preliminary test, Variant 2 provided better results. For Variant III, the $k$-means and the $k$-median clusters were solved using \textit{The C clustering library}, see \cite{Ccluster} for further details. We tested performance  using  the $k$-means and the $k$-median clustering algorithms, establishing a different number of clusters. Specifically, we carried out experiments by defining the subsets of each class as 5\%, 10\%, 20\%, 30\%, 40\%, and 50\% of the number of individuals of each class. For example, in a dataset with 100 individuals where the ratio of elements in each class is 60/40, the 10\% of clusters means making 6 clusters and 4 clusters for each class respectively. 

Table \ref{tab:StIIVariantsIjcnn} describes the best case performance (from among the different sizes of clusters and the clustering algorithms), i.e., when the problem is solved to optimality, the percentage that took the least time and the applied clustering algorithm are shown, otherwise, the ones providing the best GAP are selected. Therefore, the structure of  Table \ref{tab:StIIVariantsIjcnn} is quite similar to the previous one: the first column states the name of the dataset and the second column shows the value of parameter $C$. The following eight columns contain information about the strategy performance and the resolution process of Variants I and II. The following six columns show information about the strategy performance of Variant III: with the first stating the cluster algorithm used; the next showing the percentage chosen to determine the number of clusters; the third showing the improvement of M during the procedure; the fourth showing the strategy time plus the time spent in clustering the data; the fifth showing the total time, i.e., the clustering time, the strategy time, and the resolution time; and the sixth gives the MIP relative GAP within the time limit. The next groups of columns provide information about the behavior of Ind. Const. (LIC). Two blocks of columns provide information about this approach: in the first, the time limit is two hours and in the second, the time limit is the maximum between two hours and total time spent on the best performing Variant of Algorithm 2.2.  When the problem is solved to optimality, the procedure that takes the least time is highlighted, otherwise, the procedure that provides the best GAP within the time limit is shown. Note that we did not established a time limit for the strategies.

\begin{sidewaystable}
	\resizebox{\linewidth}{!}{
\begin{tabular}{|c|r|rrrr|rrrr|lrrrrr|rr|rr|}
	\hline
		\multicolumn{1}{|l|}{\multirow{2}[2]{*}{Data}} & \multicolumn{1}{r|}{\multirow{2}[2]{*}{$C$}} & \multicolumn{4}{c|}{Algorithm \ref{a:unbounded}.2.I} & \multicolumn{4}{c|}{Algorithm \ref{a:unbounded}.2.II} & \multicolumn{6}{c|}{Algorithm \ref{a:unbounded}.2.III} & \multicolumn{2}{c|}{Ind. Const. (LIC)}& \multicolumn{2}{c|}{Ind. Const. (LIC)}\\
			\multicolumn{1}{|l|}{} & \multicolumn{1}{r|}{} & \multicolumn{4}{c|}{} & \multicolumn{4}{c|}{} & \multicolumn{6}{c|}{} & \multicolumn{2}{c|}{T.L.:7200 s}& \multicolumn{2}{c|}{{T.L.: Best variant}}\\
				\multicolumn{1}{|l|}{} & \multicolumn{1}{r|}{} & \multicolumn{4}{c|}{} & \multicolumn{4}{c|}{} & \multicolumn{6}{c|}{} & \multicolumn{2}{c|}{}& \multicolumn{2}{c|}{{of Alg. 2}}\\
     &       & \multicolumn{1}{r}{M's Impr. } & \multicolumn{1}{r}{$t_{st}$} & \multicolumn{1}{r}{$t_{total}$} & \multicolumn{1}{r|}{GAP} & \multicolumn{1}{r}{M's Impr. } & \multicolumn{1}{r}{$t_{st}$} & \multicolumn{1}{r}{$t_{total}$} & \multicolumn{1}{r|}{GAP} & Cluster & \multicolumn{1}{r}{Per} & \multicolumn{1}{r}{M's Impr. } & \multicolumn{1}{r}{$t_{st}+t_{cl}$} & \multicolumn{1}{r}{$t_{total}$} & \multicolumn{1}{r|}{GAP} & \multicolumn{1}{r}{$t$} & \multicolumn{1}{r|}{GAP} & \multicolumn{1}{r}{$t$} & \multicolumn{1}{r|}{GAP}\\
	\hline
	\multirow{5}[2]{*}{\begin{turn}{90}Ijcnn1\_2000\end{turn}}  & 100   & 110.76\% & 12.65 & 12.78 & 0.00\% & 75.55\% & 0.34  & 0.56  & 0.00\% &   $k$-means      &   0.05    &   94.94\%
	   &  2.56
	     &    2.76
	   &  0.00\%
	     & \textbf{0.17} & \textbf{0.00\%}& \textbf{0.17} & \textbf{0.00\%} \\
 & 10    & 110.75\% & 70.64 & 70.75 & 0.00\% & 57.60\% & 0.48  & 1.33  & 0.00\% &    $k$-means    &   0.05  &  93.52\% &    4.84 &  5.11 &  0.00\%
	     & \textbf{0.39} & \textbf{0.00\%}& \textbf{0.39} & \textbf{0.00\%} \\
 & 1     & 106.07\% & 231.43 & 231.62 & 0.00\% & \textbf{42.04\%} & \textbf{0.53} & \textbf{3.54} & \textbf{0.00\%} & $k$-means & 0.05  & 87.52\% & 7.51  & 8.64 & 0.00\% & 86.90 & 0.00\% & 86.90 & 0.00\%\\
	& 0.1   & \textbf{92.90\%} & \textbf{246.30} & \textbf{618.37} & \textbf{0.00\%} & 81.48\% & 2.02  & 7214.31 & 6.05\% &  $k$-means     &    0.2   &   92.91\%
	    &   48.79    &   7255.61
	    &    4.61\%
	       & 7202.15 & 71.84\% & 7202.15 & 71.84\%\\
 & 0.01  & 75.08\% & 511.19 & 7714.37 & 12.65\% & 31.35\% & 1.33  & 7204.47 & 33.22\% &   
	  \textbf{$\pmb{k}$-median}  &  \textbf{0.4}   &  \textbf{74.50\%}
	       &    \textbf{363.30}
	       &   \textbf{7577.83}   &   \textbf{7.65\%}
	           & 7202.85 & 93.88\% & 7580.57 & 93.88\%\\

		\hline
 \multirow{5}[2]{*}{\begin{turn}{90}Ijcnn1\_3500\end{turn}}  & 100   & 99.83\% & 3015.88 & 3017.33 & 0.00\% & 11.89\% & 1.03  & 4.41  & 0.00\% & $k$-means & 0.05  & 71.46\% & 17.42 & 20.87 & \multicolumn{1}{r|}{0.00\%} & \textbf{0.79} & \textbf{0.00\%} & \textbf{0.79} & \textbf{0.00\%} \\
& 10    & 100.45\% & 2463.74 & 2464.90 & 0.00\% & 17.96\% & 1.14  & 34.97 & 0.00\% & \textbf{$\pmb{k}$-means} & \textbf{0.05} & \textbf{76.45\%} & \textbf{21.99} & \textbf{32.94} & \multicolumn{1}{r|}{\textbf{0.00\%}} & 90.64 & 0.00\% & 90.64 & 0.00\% \\
 & 1     & 99.72\% & 3834.23 & 3863.66 & 0.00\% & 23.26\% & 1.33  & 7208.13 & 20.68\% & \textbf{$\pmb{k}$-means} & \textbf{0.4} & \textbf{95.09\%} & \textbf{958.97} & \textbf{1185.13} & \multicolumn{1}{r|}{\textbf{0.00\%}} & 7207.04 & 47.34\% & 7207.04 & 47.34\% \\
 & 0.1   & 105.72\% & 4204.94 & 4205.96 & 0.00\% & 35.93\% & 2.64  & 7220.43 & 38.41\% & \textbf{$\pmb{k}$-medians} & \textbf{0.2} & \textbf{97.85\%} & \textbf{582.63} & \textbf{602.70} & \textbf{0.00\%} & 7202.48 & 83.45\% & 7202.48 & 83.45\% \\

 & 0.01  & 79.82\% & 2853.94 & 10059.62 & 27.80\% & 31.64\% & 6.16  & 7209.88 & 36.16\% & \textbf{$\pmb{k}$-median} & \textbf{0.2} & \textbf{76.49\%} & \textbf{678.20} & \textbf{7887.54} & \textbf{24.09\%} & 7202.71 & 96.04\% & 7890.85& 96.04\%\\
	\hline
\multirow{5}[2]{*}{\begin{turn}{90}Ijcnn1\_5000\end{turn}} & 100   & 99.76\% & 9715.88 & 9721.47 & 0.00\% & 10.96\% & 1.86  & 198.49 & 0.00\% & \textbf{$\pmb{k}$-means} & \textbf{0.05} & \textbf{75.10\%} & \textbf{53.85} & \textbf{165.84} & \textbf{0.00\%} & 406.17 & 0.00\% & 406.17 & 0.00\% \\
 & 10    & 99.73\% & 9228.51 & 9285.11 & 0.00\% & 13.95\% & 1.97  & 5482.24 & 0.00\% & \textbf{$\pmb{k}$-means} & \textbf{0.05} & \textbf{77.13\%} & \textbf{79.84} & \textbf{799.48} & \textbf{0.00\%} & 7202.35 & 37.50\% & 7202.35 & 37.50\% \\
 & 1     & \textbf{99.50\%} & \textbf{11661.60} & \textbf{18867.71} & \textbf{7.77\%} & 19.74\% & 2.56  & 7214.08 & 46.08\% & $k$-means & 0.5   & 95.97\% & 3547.25 & 10754.09 & 16.43\% & 7205.38 & 65.69\% &   18881.80
    &  63.32\%
\\
 & 0.1   & 104.69\% & 11977.60 & 11978.90 & 0.00\% & 32.72\% & 3.03  & 7216.60 & 46.93\% & \textbf{$\pmb{k}$-medians} & \textbf{0.5} & \textbf{100.86\%} & \textbf{4268.44} & \textbf{4315.35} & \textbf{0.00\%} & 7202.62 & 87.72\% & 7202.62 & 87.72\% \\

 & 0.01  & 78.71\% & 9264.50 & 16471.30 & 32.38\% & 32.21\% & 10.90 & 7213.46 & 50.00\% & \textbf{$\pmb{k}$-means} & \textbf{0.4} & \textbf{77.47\%} & \textbf{3810.58} & \textbf{11017.05} & \textbf{32.17\%} & 7203.10 & 96.68\%& 11021.00 & 96.68\%\\
	\hline
\end{tabular}%
}
\caption{Performance of exact approaches for solving (RL-$\ell_1$-M) and (RL-$\ell_1$) on Ijcnn.}
\label{tab:StIIVariantsIjcnn}%
\end{sidewaystable}

As shown in Table \ref{tab:StIIVariantsIjcnn}, Variant III allows us to find an equilibrium  between the time taken by the strategy and the resolution time. Variant III is especially significant in datasets with a large number of individuals, seeing as the strategy time of Variant I with the same set is huge. In the case of Ijcnn1\_5000 when parameter $C$ is equal to $0.1,$ the time taken by Variant I strategy is around three and a half hours, but when using Variant III, the strategy takes less than one hour and a half. Furthermore, the problem is solved to optimality in less than one minute after applying the strategy, while the Ind. Const. (LIC) approach provides a GAP of 87.72\% in two hours. In almost all cases, Variant III provides the best results for these datasets and when Variant I finds the optimal solution within the time limit, Variant III tends to find it in less time.  On the other hand, we noted that Variant I provides better bounds on the big M parameters than Variants II and III (the improvement of M is greater in all cases). However, the time it takes means that it is impractical. 
Whatsmore, even for a small number of clusters, i.e., 5\%, Variant III provides a significant improvement of the big M parameters compared to Variant II. 
In conclusion, we have proposed several strategies that allow us to solve many more problems to optimality within the time limit than if the indicator constraints had been used, i.e., by applying the Ind. Const. (LIC) approach. In the large dataset in which there is a high computational cost for the proposed resolution methods, we compared three different variants of Algorithm \ref{a:unbounded}. The objective of this was to achieve an equilibrium between the time employed in the strategy (which is directly related to the improvement of the big M parameters) and the time spent in solving the problem.

\subsection{Exact procedure for (RL-$\pmb{\ell_2}$-M)}

In this section, the different variants of the exact approaches proposed in Section \ref{sec:l2norm} for solving the \text{(RL-$\ell_2$-M)} formulation are tested and compared with the (RL-$\ell_2$) formulation. Specifically, we compare performance using several simulated and real-life datasets with three different resolution methods: a) the variants of Algorithm~\ref{algo:l2} -- we analyze the improvement of the big M parameters when applying the strategy, similarly to the $\ell_1$-norm case; b) the resolution method proposed in \cite{bonami} taking advantage of our strategies to tighten the values of big M parameters (since this procedure needs an unique initial big M parameter, we establish it as $\displaystyle \mbox{M}=\max_{i\in N}\{\mbox{M}_i\}$, where M$_i$ are the best values obtained with our strategy, Algorithm \ref{algo:l2}.I); c) the resolution method proposed in \cite{bonami} using Corollary \ref{coro:initialMl2} to establish the initial big M parameter; d) solving the (RL-$\ell_2$) formulation by applying the method developed in \cite{bonami} which is included in IBM-Cplex 12.7.0, i.e., generating locally valid implied bound cuts. This  latter solution approach is denoted as Ind. Const. (LIC). 

We first compare the four different methods for solving  \text{(RL-$\ell_2$-M)} and (RL-$\ell_2$) on simulated data, using Variant I of Algorithm~\ref{algo:l2}. The results are shown in Table \ref{tab: StIIInorm2}, with the first column indicating the name of the dataset and the second column the values of parameter $C$. The following groups of columns contain information about the strategy's performance and the resolution process:  a) the improvement of M during the procedure; b) the time taken by the strategy; c) the total time, i.e., strategy time plus resolution time and; d) the MIP relative GAP within the time limit (7200 seconds).
The next group of columns presents information about the resolution method proposed in \cite{bonami} but using the big M parameter obtained after applying our strategy: a) the percentage of improvement of the big M parameter with respect to the one obtained by applying Corollary \ref{coro:initialMl2}, b) the time taken by the strategies; c) the total time, i.e., strategy time plus resolution time and; d) the MIP relative GAP within the time limit. 
The third group of columns report the results of the same method of \cite{bonami}  but with the initial big M parameter (Corollary \ref{coro:initialMl2}): a) the time taken by the strategy; b) the total time, i.e., strategy time plus resolution time and; c) the MIP relative GAP within the time limit.
 The last group of columns shows the behavior of Ind. Const. (LIC): the first column shows the total time and the second shows the GAP obtained  within the time limit. The strategy with the best performance is highlighted.
 
\begin{table}[htbp]
	\centering
	\resizebox{\hsize}{!}{		\begin{tabular}{|c|r|rrrr|rrrr|rrr|rr|}
			\hline
			\multicolumn{1}{|l|}{\multirow{2}[2]{*}{Data}} & \multicolumn{1}{r|}{\multirow{2}[2]{*}{$C$}} & \multicolumn{4}{c|}{Algorithm \ref{algo:l2}.I}& \multicolumn{4}{c|}{\cite{bonami}} & \multicolumn{3}{c|}{\cite{bonami} } & \multicolumn{2}{c|}{Ind. Const. (LIC)} \\
			&	&  &&&   & \multicolumn{4}{c|}{Big M: Algorithm \ref{algo:l2}.I} & \multicolumn{3}{c|}{Big M: Corollary \ref{coro:initialMl2}} &&\\
			&       & \multicolumn{1}{r}{M's Impr. } & \multicolumn{1}{r}{$t_{st}$} & \multicolumn{1}{r}{$t_{total}$} & \multicolumn{1}{r|}{GAP} &\multicolumn{1}{r}{M's Impr. }&\multicolumn{1}{r}{$t_{st}$} & \multicolumn{1}{r}{$t_{total}$}& \multicolumn{1}{r|}{GAP} & \multicolumn{1}{r}{$t_{st}$} & \multicolumn{1}{r}{$t_{total}$}& \multicolumn{1}{r|}{GAP} & \multicolumn{1}{r}{$t$} & \multicolumn{1}{r|}{GAP} \\
			\hline
\multirow{5}[2]{*}{\begin{turn}{90} n160d2A\end{turn}}  & 100   & \textbf{95.13\%} & \textbf{21.20} & \textbf{23.28} & \textbf{0.00\%} & 58.31\% & 30.43  & 7230.49 & 52.94\% & 9.80  & 7209.83 & 54.38\% & 7204.46 & 57.71\% \\
 & 10    & \textbf{95.98\%} & \textbf{20.21} & \textbf{20.47} & \textbf{0.00\%}     & 68.09\%& 30.21 &473.61 & 0.00\% & 7.63  & 7207.74 & 58.49\% & 7203.04 & 69.05\% \\
 & 1     & \textbf{92.80\%} & \textbf{33.82} & \textbf{33.86} & \textbf{0.00\%}     & 75.46\%& 46.92 &7247.01 & 25.02\%& 12.34 & 7212.40 & 59.52\% & 7202.29 & 79.63\% \\
 & 0.1   & \textbf{34.31\%} & \textbf{29.94} & \textbf{7232.15} & \textbf{35.63\%}  & 37.01\%& 40.88 &7240.91 & 59.67\%& 9.75  & 7209.80 & 65.12\% & 7202.11 & 83.21\% \\
 & 0.01  & \textbf{37.43\%} & \textbf{25.06} & \textbf{7227.78} & \textbf{0.87\%}   & 38.91\%& 36.19 &7236.22 & 3.27\% & 8.84  & 7208.92 & 13.47\% & 7202.82 & 84.39\% \\
\hline
\multirow{5}[2]{*}{\begin{turn}{90} n160d5A\end{turn}} & 100   & \textbf{91.82\%} & \textbf{21.92} & \textbf{7231.52} & \textbf{19.89\%} & 52.79\% & 33.76 & 7233.81 & 60.78\% & 10.24 & 7210.29 & 66.21\% & 7202.70 & 60.16\% \\
 & 10    & \textbf{92.63\%} & \textbf{25.31} & \textbf{920.06} & \textbf{0.00\%}  & 60.55\% &40.12 & 7214.93 & 53.72\%& 11.17 & 7240.24 & 68.06\% & 7202.57 & 71.49\% \\
 & 0.1   & \textbf{40.67\%} & \textbf{34.60} & \textbf{7236.97} & \textbf{41.55\%}& 60.46\% &43.96 & 7209.43 & 46.69\%& 9.39  & 7244.03 & 52.77\% & 7201.94 & 84.20\% \\
 & 1     & \textbf{87.19\%} & \textbf{38.24} & \textbf{7238.28} & \textbf{2.00\%} & 37.42\% &46.75 & 7208.65 & 29.77\%& 10.37 & 7246.89 & 65.56\% & 7202.43 & 79.69\% \\
 & 0.01  & \textbf{36.06\%} & \textbf{40.81} & \textbf{7242.89} & \textbf{3.16\%} & 51.97\% &53.65 & 7212.95 & 8.28\% & 12.14 & 7253.76 & 52.44\% & 7202.49 & 85.67\% \\
\hline
\multirow{5}[2]{*}{\begin{turn}{90} n160d10A\end{turn}} & 100   & \textbf{88.80\%} & \textbf{23.80} & \textbf{7225.69} & \textbf{50.92\%}  & 47.88\% &  31.99 & 7232.31 & 81.70\% & 7.51  & 7207.62 & 84.98\% & 7202.79 & 69.39\% \\
 & 10    & \textbf{88.93\%} & \textbf{26.71} & \textbf{7228.74} & \textbf{39.06\%} & 51.62\% & 37.27& 7210.62 & 71.99\% & 8.27  & 7237.33 & 76.19\% & 7203.00 & 73.20\% \\
 & 1     & \textbf{80.82\%} & \textbf{44.99} & \textbf{7246.83} & \textbf{49.48\%} & 45.03\% & 54.31& 7209.33 & 69.49\% & 7.67  & 7254.32 & 72.83\% & 7205.35 & 82.99\% \\
 & 0.1   & \textbf{32.85\%} & \textbf{44.83} & \textbf{7246.93} & \textbf{55.43\%} & 23.53\% & 58.27& 7215.04 & 65.28\% & 9.79  & 7259.87 & 65.23\% & 7204.11 & 84.65\% \\
 & 0.01  & \textbf{22.29\%} & \textbf{42.85} & \textbf{7244.36} & \textbf{27.69\%} & 18.80\% & 53.88& 7211.13 & 29.37\% & 11.69 & 7253.98 & 35.19\% & 7206.03 & 85.88\% \\
\hline
\multirow{5}[2]{*}{\begin{turn}{90} n160d2B\end{turn}} & 100   & \textbf{45.72\%} & \textbf{34.55} & \textbf{7234.60} & \textbf{14.72\%}  & 37.03\% & 43.94  & 7243.96 & 65.37\% & 8.96  & 7209.04 & 67.99\% & 7207.06 & 63.91\% \\
 & 10    & \textbf{45.71\%} & \textbf{32.08} & \textbf{7232.14} & \textbf{7.86\%}& 37.11\% &42.43 &  7242.54 & 62.53\% & 9.67  & 7209.68 & 63.03\% & 7211.17 & 69.15\% \\
 & 1     & \textbf{45.83\%} & \textbf{27.53} & \textbf{6342.20} & \textbf{0.00\%}& 37.56\% &36.90 &  7236.91 & 5.03\%  & 7.94  & 7208.03 & 55.79\% & 7206.39 & 74.79\% \\
 & 0.1   & 49.92\%          & 25.93          & 7226.06          & 8.71\%         & 46.79\% &34.96 &  7234.98 & 15.03\% & \textbf{7.65} & \textbf{7207.77}  & \textbf{6.95\%} & 7202.59 & 82.60\% \\
 & 0.01  & \textbf{99.95\%} & \textbf{12.53} & \textbf{12.55} & \textbf{0.00\%}  & \textbf{99.98\%} &\textbf{12.55} & \textbf{12.55}  & \textbf{0.00\%}    & 10.83 & 7210.89 & 8.47\% & 7207.95 & 82.65\% \\
\hline
\multirow{5}[2]{*}{\begin{turn}{90} n160d5B\end{turn}} & 100   & \textbf{51.88\%} & \textbf{36.08} & \textbf{7236.26} & \textbf{63.96\%}  & 34.82\% &  44.90 & 7244.95 & 72.96\%& 7.93  & 7207.97 & 73.86\% & 7207.12 & 65.73\% \\
 & 10    & \textbf{51.50\%} & \textbf{32.79} & \textbf{7232.90} & \textbf{50.58\%}& 34.81\% & 44.22  & 7244.23 & 74.30\% & 11.18 & 7211.23 & 71.17\% & 7204.61 & 72.91\% \\
 & 1     & \textbf{44.65\%} & \textbf{29.60} & \textbf{7229.63} & \textbf{36.32\%}& 33.80\% & 38.59  & 7238.61 & 52.32\% & 10.96 & 7211.09 & 49.82\% & 7204.81 & 75.30\% \\
 & 0.1   & \textbf{31.88\%} & \textbf{23.45} & \textbf{7225.80} & \textbf{47.80\%}& 31.55\% & 31.56  & 7232.10 & 41.38\% & 16.17 & 7216.26 & 54.28\% & 7206.88 & 81.36\% \\
 & 0.01  & \textbf{50.31\%} & \textbf{33.46} & \textbf{7236.20} & \textbf{6.06\%} & 57.05\% & 40.83  & 7240.89 & 17.02\% & 7.28  & 7207.39 & 53.45\% & 7202.51 & 83.54\% \\
\hline
\multirow{5}[2]{*}{\begin{turn}{90} n160d10B\end{turn}} & 100   & \textbf{30.36\%} & \textbf{28.84} & \textbf{7234.94} & \textbf{51.97\%} & 32.10\% & 44.32  & 7244.52 & 80.00\%& 16.31 & 7216.37 & 82.20\% & 7202.03 & 65.16\% \\
 & 10    & \textbf{29.70\%} & \textbf{32.81} & \textbf{7234.46} & \textbf{53.72\%} & 31.95\% & 43.32  & 7243.42 & 78.16\% & 8.31  & 7208.51 & 79.61\% & 7217.02 & 68.35\% \\
 & 1     & \textbf{27.33\%} & \textbf{28.70} & \textbf{7230.35} & \textbf{51.52\%} & 31.04\% & 45.42  & 7245.52 & 67.32\% & 7.23  & 7207.26 & 72.71\% & 7202.42 & 70.35\% \\
 & 0.1   & \textbf{25.35\%} & \textbf{29.57} & \textbf{7231.70} & \textbf{50.20\%} & 30.23\% & 45.23  & 7245.44 & 60.36\% & 12.55 & 7214.25 & 70.18\% & 7202.04 & 77.84\% \\
 & 0.01  & \textbf{25.02\%} & \textbf{25.67} & \textbf{7229.13} & \textbf{33.20\%} & 32.80\% & 36.14  & 7236.21 & 49.54\% & 14.77 & 7214.86 & 60.16\% & 7202.42 & 84.91\% \\
\hline
		\end{tabular}%
	}
	\caption{Performance of exact approaches on simulated data for solving (RL-$\ell_2$-M) and (RL-$\ell_2$).}
	\label{tab: StIIInorm2}%
\end{table}

Note that, in almost all the cases, the behavior of our strategy is better than  the algorithm proposed in \cite{bonami}, which in turn is better than solving the problem using the Ind. Const. (LIC) approach. For instance, in dataset n160d2A with a $C$ parameter equal to 1, our strategy provides the optimal solution in a total time of 34 seconds.
Meanwhile, the GAP of the algorithm from \cite{bonami} using the improved big M derived from Algorithm \ref{algo:l2}.I and with a time limit of two hours is 25.02\%, and the same algorithm without the improved big M is 59.02\%. Finally, the GAP of the Ind. Const. (LIC) approach is 79.63\%. Observe that the strategies take less than one minute in all cases. It should also be remarked that the improvement of algorithm of \cite{bonami} provides better final GAPs when Algorithm \ref{algo:l2}.I is used to initialize the big M parameter.

Next, we compare the performance of the three strategies using real-life datasets. The results are reported in Table \ref{tab: StIIInorm2_sim} which has the same structure as Table \ref{tab: StIIInorm2}. 

\begin{table}[htbp]
	\centering
	\resizebox{\hsize}{!}{
\begin{tabular}{|c|r|rrrr|rrrr|rrr|rr|}
	\hline
	\multicolumn{1}{|l|}{\multirow{2}[2]{*}{Data}} & \multicolumn{1}{r|}{\multirow{2}[2]{*}{$C$}} & \multicolumn{4}{c|}{Algorithm \ref{algo:l2}.I}& \multicolumn{4}{c|}{\cite{bonami}} & \multicolumn{3}{c|}{\cite{bonami} } & \multicolumn{2}{c|}{Ind. Const. (LIC)} \\
		&	&  &&&   & \multicolumn{4}{c|}{Big M: Algorithm \ref{algo:l2}.I} & \multicolumn{3}{c|}{Big M: Corollary \ref{coro:initialMl2}} &&\\
	&       & \multicolumn{1}{r}{M's Impr. } & \multicolumn{1}{r}{$t_{st}$} & \multicolumn{1}{r}{$t_{total}$} & \multicolumn{1}{r|}{GAP} &\multicolumn{1}{r}{M's Impr. }&\multicolumn{1}{r}{$t_{st}$} & \multicolumn{1}{r}{$t_{total}$}& \multicolumn{1}{r|}{GAP} & \multicolumn{1}{r}{$t_{st}$} & \multicolumn{1}{r}{$t_{total}$}& \multicolumn{1}{r|}{GAP} & \multicolumn{1}{r}{$t$} & \multicolumn{1}{r|}{GAP} \\
	\hline
\multirow{5}[2]{*}{\begin{turn}{90} Wpbc\end{turn}}   & 100   & \textbf{75.03\%} & \textbf{83.20} & \textbf{7283.31} & \textbf{61.50\%}  & 57.31\% &113.60 & 7313.72 & 89.49\% & 29.02 & 7229.30 & 91.18\% & 7200.30 & 74.89\% \\
 & 10    & \textbf{73.12\%} & \textbf{91.57} & \textbf{7298.90} & \textbf{54.19\%}& 53.52\% &137.31 & 7246.03 & 85.81\%& 31.97 & 7232.50 & 88.93\% & 7337.60 & 79.39\% \\
 & 1     & \textbf{62.64\%} & \textbf{78.19} & \textbf{7280.04} & \textbf{56.31\%}& 35.68\% &121.74 & 7243.93 & 81.73\%& 30.26 & 7231.83 & 84.48\% & 7322.12 & 79.00\% \\
 & 0.1   & \textbf{61.73\%} & \textbf{64.19} & \textbf{7265.96} & \textbf{54.23\%}& 31.12\% &119.18 & 7255.07 & 71.00\%& 49.84 & 7250.18 & 73.79\% & 7319.26 & 78.06\% \\
 & 0.01  & \textbf{63.82\%} & \textbf{53.60} & \textbf{7255.50} & \textbf{37.44\%}& 30.47\% &84.93  & 7231.53 & 46.73\%& 30.76 & 7230.86 & 55.59\% & 7285.13 & 76.30\% \\
\hline
\multirow{5}[2]{*}{\begin{turn}{90} SONAR\end{turn}}  & 100   & \textbf{90.09\%} & \textbf{162.79} & \textbf{173.73} & \textbf{0.00\%}& 79.77\% &392.81 & 7593.58 & 44.02\%  & 325.01 & 7525.11 & 55.63\% & 511.82 & 0.00\% \\
 & 10    & \textbf{78.74\%} & \textbf{149.84} & \textbf{7350.08} & \textbf{17.15\%}& 64.47\% &286.62 & 7337.49 & 76.65\% & 223.61 & 7423.68 & 82.16\% &7487.33& 62.90\% \\
 & 1     & \textbf{73.85\%} & \textbf{119.74} & \textbf{7321.56} & \textbf{39.72\%}& 58.61\% &319.97 & 7400.57 & 75.73\% & 229.65 & 7429.91 & 83.53\% &7520.31& 75.27\% \\
 & 0.1   & \textbf{64.39\%} & \textbf{115.54} & \textbf{7317.44} & \textbf{50.75\%}& 49.44\% &197.02& 7282.00 & 74.16\%  & 70.71 & 7270.97 & 82.77\% & 7397.54& 82.49\% \\
 & 0.01  & \textbf{50.34\%} & \textbf{147.20} & \textbf{7348.83} & \textbf{47.62\%}& 36.81\% &206.31& 7259.84 & 66.96\%  & 56.31 & 7257.78 & 75.84\% & 7407.04& 88.27\% \\
\hline
\multirow{5}[2]{*}{\begin{turn}{90} SPECT\end{turn}}  & 100   & \textbf{70.21\%} & \textbf{189.01} & \textbf{836.86} & \textbf{0.00\%}& 38.47\% &203.93 & 7404.06 & 24.99\%  & 16.10 & 7216.29 & 29.13\% & 7200.23 & 76.40\% \\
 & 10    & 70.13\% & 158.47 & 7361.31 & 35.42\% & \textbf{38.50\%} &\textbf{183.31} & \textbf{7383.41} & \textbf{33.16\%}&  22.40 & 7222.62 & 35.03\%&7200.29 & 78.39\% \\
 & 1     & 69.89\% & 150.57 & 7353.08 & 51.00\% & \textbf{38.56\%} &\textbf{183.86} & \textbf{7383.96} & \textbf{33.91\%}&  26.49 & 7226.53 &38.57\%&7200.35 & 79.05\% \\
 & 0.1   & 65.31\% & 170.15 & 7376.85 & 36.75\%                  & \textbf{34.11\%} &\textbf{204.98} & \textbf{7405.02} & \textbf{31.65\%}& 21.31 & 7221.50 & 42.62\% & 7200.19 & 79.81\% \\
 & 0.01  & 55.10\% & 95.92 & 7298.00 & 53.97\% & \textbf{14.43\%} &\textbf{113.45} & \textbf{7313.56} & \textbf{41.12\%}& 14.01 & 7214.08 & 41.58\% & 7209.16 & 80.60\%\\
\hline
\multirow{5}[2]{*}{\begin{turn}{90} IONO\end{turn}}   & 100   & \textbf{88.33\%} & \textbf{266.56} & \textbf{7467.05} & \textbf{22.43\%} & 65.76\% & 398.20 & 7598.75 & 31.49\% & 118.78 & 7319.17 & 33.80\% & 7200.40 & 58.52\% \\
  & 10    & \textbf{85.44\%} & \textbf{245.85} & \textbf{7446.34} & \textbf{35.13\%} & 61.24\% & 398.64 &7599.30 & 39.67\%& 209.33 & 7409.57 & 42.83\% & 7200.96 & 65.65\% \\
  & 1     & \textbf{81.51\%} & \textbf{228.17} & \textbf{7429.80} & \textbf{42.81\%} & 57.00\% & 306.68& 7507.11& 43.37\% & 92.62 & 7293.17 & 49.11\% & 7200.51 & 78.51\% \\
  & 0.1   & \textbf{78.84\%} & \textbf{207.92} & \textbf{7411.75} & \textbf{40.79\%} & 56.07\% & 255.25& 7460.68& 47.12\% & 48.75 & 7250.88 & 53.06\% & 7201.60 & 83.29\% \\
  & 0.01  & \textbf{69.17\%} & \textbf{214.34} & \textbf{7430.09} & \textbf{48.28\%} & 45.62\% & 258.26& 7458.58& 50.96\% & 40.51 & 7240.53 & 58.76\% & 7200.25 & 90.34\% \\
\hline
\multirow{5}[2]{*}{\begin{turn}{90} Wdbc\end{turn}}   & 100   & 97.51\% & 227.89 & 230.00 & 0.00\% & 81.18\% & 289.28 & 289.28 & 0.00\% & \textbf{112.28} & \textbf{112.28} & \textbf{0.00\%} & 249.77 & 0.00\% \\
  & 10    & \textbf{98.62\%} & \textbf{319.62} & \textbf{325.70} & \textbf{0.00\%} & 83.09\% & 412.12 & 7612.31 & 27.02\% & 207.31 & 7407.50 & 43.83\% & 7200.44 & 34.59\% \\
  & 1     & \textbf{103.31\%} & \textbf{398.20} & \textbf{398.62} & \textbf{0.00\%}& 87.79\% & 439.49 & 439.49  & 0.00\%   & 93.88 & 7296.41 & 58.05\% & 7200.28 & 68.70\% \\
  & 0.1   & \textbf{103.59\%} & \textbf{410.37} & \textbf{410.53} & \textbf{0.00\%}& 83.46\% & 421.18 & 421.18  & 0.00\%   & 61.46 & 7262.07 & 56.14\% & 7200.49 & 87.19\% \\
  & 0.01  & \textbf{101.29\%} & \textbf{403.58} & \textbf{403.63} & \textbf{0.00\%}& 71.79\% & 416.01 & 416.01  & 0.00\%   & 61.18 & 7261.19 & 57.22\% & 7200.27 & 95.10\% \\
\hline
\multirow{5}[2]{*}{\begin{turn}{90} WBC\end{turn}}   & 100   & \textbf{93.37\%} & \textbf{244.73} & \textbf{313.64} & \textbf{0.00\%}& 68.10\% & 296.78 & 7497.08 & 34.13\% & 46.19 & 7246.58 & 51.41\% & 7200.33 & 40.99\% \\
   & 10    & \textbf{93.65\%} & \textbf{211.91} & \textbf{271.58} & \textbf{0.00\%} & 68.34\% & 253.61 & 7453.80 & 29.22\% & 53.19 & 7253.63 & 43.74\% & 7200.50 & 43.27\% \\
   & 1     & \textbf{94.49\%} & \textbf{214.41} & \textbf{277.75} & \textbf{0.00\%} & 68.99\% & 249.87 & 7450.14 & 22.38\% & 51.34 & 7251.51 & 38.64\% & 7200.71 & 43.56\% \\
   & 0.1   & \textbf{99.24\%} & \textbf{355.86} & \textbf{356.35} & \textbf{0.00\%} & 72.64\% & 360.38 & 360.37  & 0.00\%    & 30.73 & 7230.96 & 27.19\% & 7200.31 & 60.36\% \\
   & 0.01  & 105.06\% & 364.61 & 364.69 & 0.00\%                                    & 70.68\% & 369.28 & 369.28  & 0.00\%    & \textbf{9.47} & \textbf{9.47} & \textbf{0.00\%} & 7200.27 & 80.36\% \\
	\hline
\end{tabular}%

}
\caption{Performance of exact approaches on real data for solving (RL-$\ell_2$-M) and (RL-$\ell_2$).}
\label{tab: StIIInorm2_sim}%
\end{table}

Table \ref{tab: StIIInorm2_sim} shows that behavior using real-life datasets is quite similar to simulated datasets, i.e., our strategy performs better in almost all the cases. The   Wdbc and WBC datasets are of particular interest as the problem is solved in less than 411 seconds in all the cases. However, four cases remain unsolved after two hours when using the method proposed in \cite{bonami} using the improved big M parameters derived from Algorithm \ref{algo:l2}.I, providing GAPs of between 22.38\% and 34.13\%. In addition, eight cases remain unsolved after two hours when using the method in \cite{bonami} using Corollary \ref{coro:initialMl2} and nine cases when using the Ind. Const. (LIC). These two last methods provide GAPs of between  27.19\% and 95.10\%.

Lastly, we analyze three large datasets: Ijcnn1\_2000, Ijcnn1\_3500, and Ijcnn1\_5000. We compare the performance of Variants I and II of Algorithm \ref{algo:l2} with the performance of Ind. Const. (LIC) and with the resolution method proposed in \cite{bonami}. The results of Variant III have been omitted because it performed worse than Variant II in all the cases, i.e., it did not improve the GAP or the overall time. 

The results are shown in Table \ref{tb:Al3large}, which has a similar structure to Table \ref{tab: StIIInorm2}: the first column shows the name of the dataset and the second column the value of parameter $C$. The following eleven columns contain information about the strategy's performance and the resolution process of Variant I, Variant I version two (v2) and Variant II. Version two of Variant I modifies steps 11-12 of Algorithm \ref{algo:l2}: it only solves the problems for $i\in N$ such that $M_i$ is greater than the median of the $M_i$ values. The next groups of columns show information about the behavior of the Algorithm of \cite{bonami} and the Ind. Const. (LIC) approach. 
Three blocks of columns provide information about \cite{bonami} method: the first block reports the results of this method using the improved initial big M parameters derived from our strategy; the second block is the same method using Corollary \ref{coro:initialMl2} to initialize the big M and establishing a time limit of two hours; and the last one reports the results of the method using Corollary \ref{coro:initialMl2} to initialize the big M and establishing as time limit the total time of the best performing Variant of Algorithm \ref{algo:l2} whenever said time is greater than two hours, otherwise it is two hours.  
Finally, two blocks of columns provide information about the Ind. Const. (LIC) approach: in the first, the time limit is two hours and in the second, it is the total time of the best performing Variant of Algorithm \ref{algo:l2} if it is greater than two hours or two hours otherwise. When the problem is solved to optimality, the strategy that takes the least time is shown in bold, otherwise, the strategy that provides the best GAP within the time limit is highlighted. The column for the improvement of M in Algorithm \ref{algo:l2}.II has been omitted because it is almost zero in all cases.

\begin{table}
\begin{center}

	\resizebox{\hsize}{!}{
	\begin{tabular}{|c|r|rrrr|rrrr|rrr|}
		\hline
			\multicolumn{1}{|l|}{\multirow{2}[2]{*}{Data}} & \multicolumn{1}{r|}{\multirow{2}[2]{*}{C}} & \multicolumn{4}{c|}{Algorithm \ref{algo:l2}.I} & \multicolumn{4}{c|}{Algorithm \ref{algo:l2}.I v2} &
			\multicolumn{3}{c|}{Algorithm \ref{algo:l2}.II}  \\
			&       & \multicolumn{1}{r}{M's Impr. } & \multicolumn{1}{r}{$t_{st}$} & \multicolumn{1}{r}{$t_{total}$} & \multicolumn{1}{r|}{GAP} &\multicolumn{1}{r}{M's Impr. } &  \multicolumn{1}{r}{$t_{st}$} & \multicolumn{1}{r}{$t_{total}$} & \multicolumn{1}{r|}{GAP} & \multicolumn{1}{r}{$t_{st}$} & \multicolumn{1}{r}{$t_{total}$} & \multicolumn{1}{r|}{GAP} \\
		\hline

\multirow{5}[2]{*}{\begin{turn}{90} Ijcnn1\_2000\end{turn}}  & 100   & 107.69\% & 2150.59 & 2150.61 & 0.00\% & 107.31\% & 1506.21 & 1506.23 & 0.00\% & 25.31 & 25.61 & 0.00\% \\
 & 10    & 106.74\% & 5214.58 & 5214.61 & 0.00\% & 103.04\%            & 2032.12 & 2035.20 & 0.00\% & 19.98 & 21.92 & 0.00\% \\
 & 1     & 104.42\% & 4492.10 & 4494.41 & 0.00\% & 102.40\%            & 2022.32 & 2025.60 & 0.00\% & 20.55 & 5431.52 & 0.00\% \\
 & 0.1   & 103.79\% & 5582.22 & 5583.06 & 0.00\% & \textbf{101.83\%} & \textbf{3581.91} & \textbf{3584.28} & \textbf{0.00\%} & 14.71 & 7218.11 & 55.74\% \\
 & 0.01  & 76.60\% & 2970.38 & 10175.25 & 67.46\% & \textbf{76.81\%} & \textbf{3182.42} & \textbf{10388.21} & \textbf{67.11\%} & 15.46 & 7220.78 & 73.36\% \\
\hline
\multirow{5}[2]{*}{\begin{turn}{90} Ijcnn1\_3500\end{turn}} & 100   & 97.04\% & 6990.69 & 6993.54 & 0.00\% & 49.26\% & 3617.80 & 3626.25 & 0.00\% & 69.50 & 94.29 & 0.00\% \\
 & 10    & 94.98\% & 6841.82 & 7210.65 & 0.00\% & \textbf{47.28\%} &   \textbf{2960.29} & \textbf{5446.03} & \textbf{0.00\%} & 64.77 & 7265.28 & 19.99\% \\
 & 1     & 99.98\% & 33883.80 & 33989.79 & 0.00\% & \textbf{98.81\%} & \textbf{19193.70} & \textbf{21152.66} & \textbf{0.00\%} & 63.84 & 7264.28 & 59.93\% \\
 & 0.1   & 100.73\% & 45828.20 & 45828.87 & 0.00\% & \textbf{98.92\%} &\textbf{28850.40} & \textbf{28868.85} & \textbf{0.00\%} & 61.88 & 7262.03 & 72.26\% \\
 & 0.01  & 78.69\% & 17913.20 & 25116.00 & 70.17\% & \textbf{78.63\%} &\textbf{15786.20} & \textbf{22988.27} & \textbf{70.11\%} & 45.37 & 7247.83 & 79.03\% \\
\hline
\multirow{5}[2]{*}{\begin{turn}{90} Ijcnn1\_5000\end{turn}} & 100   & 97.40\% & 17403.40 & 17799.19 & 0.00\% & \textbf{50.24\%} & \textbf{9374.54} & \textbf{10840.45} & \textbf{0.00\%} & 117.28 & 7317.57 & 9.19\% \\
 & 10    & 99.61\% & 60964.80 & 64717.57 & 0.00\% & \textbf{99.32\%} &   \textbf{44729.00} & \textbf{50071.69} & \textbf{0.00\%} & 111.71 & 7311.93 & 58.03\% \\
 & 1     & 100.13\% & 88536.30 & 95737.49 & 6.75\% & \textbf{100.19\%} & \textbf{69876.50} & \textbf{71848.85} & \textbf{0.00\%} & 136.72 & 7337.00 & 68.23\% \\
 & 0.1   & 103.16\% & 128749.00 & 128749.59 & 0.00\% & \textbf{100.85\%}&\textbf{78689.90} & \textbf{78690.51} & \textbf{0.00\%} & 117.02 & 7317.40 & 70.86\% \\
 & 0.01  & 81.79\% & 37206.80 & 44408.84 & 68.86\% & \textbf{81.38\%} &  \textbf{26749.60} & \textbf{33951.29} & \textbf{68.89\%} & 67.68 & 7267.831 & 82.22\% \\

		\hline
	\end{tabular}%
}

\bigskip
\resizebox{\hsize}{!}{
	\begin{tabular}{|c|r|rrrr|rrr|rrr|rr|rr|}
		\hline
			\multicolumn{1}{|l|}{\multirow{3}[2]{*}{Data}} & \multicolumn{1}{r|}{\multirow{3}[2]{*}{C}} & \multicolumn{4}{c|}{\cite{bonami}} &
			 \multicolumn{3}{c|}{\cite{bonami}} &
			\multicolumn{3}{c|}{\cite{bonami}} & \multicolumn{2}{c|}{Ind. Const. (LIC)}  & \multicolumn{2}{c|}{Ind. Const. (LIC)} \\
				&	&  \multicolumn{4}{c|}{Big M: Algorithm \ref{algo:l2}.I v2}   & \multicolumn{3}{c|}{Big M: Corollary \ref{coro:initialMl2}} & \multicolumn{3}{c|}{Big M: Corollary \ref{coro:initialMl2}} &\multicolumn{2}{c|}{T.L.: 7200 s } & \multicolumn{2}{c|}{T.L.: Best variant}\\
			&	&  \multicolumn{4}{c|}{T.L.: 7200 s }   & \multicolumn{3}{c|}{T.L.: 7200 s } & \multicolumn{3}{c|}{T.L.: Best variant of Alg. 3} && & \multicolumn{2}{c|}{of Alg. 3}\\
			&   &\multicolumn{1}{r}{M's Impr. } & \multicolumn{1}{r}{$t_{st}$} & \multicolumn{1}{r}{$t_{total}$} & \multicolumn{1}{r|}{GAP}& \multicolumn{1}{r}{$t_{st}$} & \multicolumn{1}{r}{$t_{total}$} & \multicolumn{1}{r|}{GAP}  & \multicolumn{1}{r}{$t_{st}$} & \multicolumn{1}{r}{$t_{total}$} & \multicolumn{1}{r|}{GAP}  & \multicolumn{1}{r}{$t_{total}$} & 
		\multicolumn{1}{r|}{GAP} & \multicolumn{1}{r}{$t_{total}$} & \multicolumn{1}{r|}{GAP}  \\

\hline
\multirow{5}[2]{*}{\begin{turn}{90} Ijcnn1\_2000\end{turn}} & 100 & 100.00\% &1506.30 & 1506.30  & 0.00\%    & 13.17 & 13.17 & 0.00\% &  13.17 & 13.17 & 0.00\% & \textbf{0.19} & \textbf{0.00\%} & \textbf{0.19} & \textbf{0.00\%}\\
 & 10   & 100.00\% &2032.20 & 2032.20  & 0.00\%  & 11.28 & 11.28 & 0.00\% &       11.28 & 11.28 & 0.00\%       & \textbf{7.61} & \textbf{0.00\%} & \textbf{7.61} & \textbf{0.00\%} \\
 & 1    & 94.66\%  &2096.32 & 2096.32  & 0.00\%  & \textbf{180.87} & \textbf{180.87} & \textbf{0.00\%} &      \textbf{180.87} & \textbf{180.87} & \textbf{0.00\%} & 7200.91 &42.98\% & 7200.91 & 42.98\% \\
 & 0.1  & 93.87\%  &3582.38 & 3585.46  & 0.00\%  & 1272.02 & 8472.07 & 60.63\% &    1272.02 & 8472.07 & 60.63\%    & 7202.76 & 85.99\% & 7202.76 & 85.99\% \\
 & 0.01 & 28.90\%  &4032.29 & 11232.78 & 67.84\% & 245.79 & 7446.41 & 79.24\% &  262.92  &  10915.31&    78.73\%   & 7202.34 & 95.29\% & 10845.10 & 95.26\% \\
\hline
\multirow{5}[2]{*}{\begin{turn}{90} Ijcnn1\_3500\end{turn}} & 100  &23.07\% &3673.52 &  3673.52 & 0.00\%   & \textbf{74.68} & \textbf{74.68} & \textbf{0.00\%} &      \textbf{74.68} & \textbf{74.68} & \textbf{0.00\%}&        432.85 & 0.00\% & 432.85 & 0.00\% \\
 & 10    &22.92\% &3998.42 & 11198.63& 54.48\% & 1068.75 & 8268.93 & 62.61\% &     1068.75 & 8268.93 & 62.61\% & 7201.19 & 53.97\% & 7201.19 & 53.97\%  \\
 & 1     &90.09\% &19649.06 & 26849.16& 29.95\% & 1657.10 & 8860.39 & 69.04\% &   1666.30   &     21161.42  &   65.34\% & 7202.72 & 78.18\% & 21156.60& 61.83\% \\
 & 0.1   &93.64\% &29455.93 & 36674.59& 7.09\%  & 1487.32 & 8687.86 & 73.70\% &   1395.12    &  28778.88     &    73.18\%  & 7201.95 & 92.10\% & 28874.30 & 91.18\% \\
 & 0.01  &36.54\% &17059.32 & 24260.11& 77.67\% & 579.66 & 7781.42 & 81.00\% &   589.97    &   23164.37    &    80.92\%   & 7202.58 & 97.43\% & 22996.00 & 97.42\% \\
\hline
\multirow{5}[2]{*}{\begin{turn}{90} Ijcnn1\_5000\end{turn}}  & 100 & 21.68\% & 11457.14 & 18657.74  & 56.98\%   & 2438.92 & 9639.52 & 62.26\% &   2456.31   & 10858.61 &    60.77\%  & 7201.14 & 57.52\% &  10842.30 & 53.73\% \\
 & 10    & 91.37\% & 45451.16&  52651.21  & 39.32\% & 3052.90 & 10270.38 & 73.14\% &    3072.17  &  50094.04     & 66.06\%  & 7201.44 & 70.90\% & 50077.30 & 65.91\% \\
 & 1     & 91.89\% & 70700.40&  77900.52  & 39.83\% & 2682.25 & 9882.31 & 76.04\% &  2692.03  & 71880.28 &69.67\% & 7201.63 & 84.37\% &   71856.40    &  81.66\% \\
 & 0.1   & 93.78\% & 79927.33 & 87129.11  & 11.96\% & 2730.07 & 9941.87 & 75.95\% &  2754.24  &   78720.34    &   72.61\%    & 7202.71 & 94.84\% & 78708.30    & 93.17\% \\
 & 0.01  & 41.07\% & 28747.69 & 35951.31  & 77.42\% & 1045.90 & 8257.49 & 83.86\% &    1051.38   &     33959.29  &  82.52\%    & 7202.60 & 97.96\% &  33958.50& 97.52\% \\
\hline
\end{tabular}%

}
    
\end{center}
\caption{Performance of exact approaches on Ijcnn for solving (RL-$\ell_2$-M) and (RL-$\ell_2$).}
\label{tb:Al3large}%
\end{table}	

Note that our strategy works well even in large datasets. In fact, Algorithm \ref{algo:l2}.I v2 outperforms Algorithm \ref{algo:l2}.I. For example, using Algorithm \ref{algo:l2}.I  v2 for the Ijcnn1\_2000 dataset, with the value of parameter $C$ being equal to 0.1 and after a preprocessing of 3581.91 seconds, we are able to solve the problem in less than three seconds. When using \cite{bonami} method, this instance can only be solved when the big M parameter is the one derived from the Algorithm  \ref{algo:l2}.I v2.
 When applying the Ind. Const. (LIC) approach, the GAP is 85.99\% after two hours. Although the time employed in the strategy is huge in some cases, we are able to solve almost all the problems to optimality. Over the same amount of time, the GAPs provided by the Algorithm of \cite{bonami} using Corollary \ref{coro:initialMl2} and the Ind. Const. (LIC) approach are considerably bigger. For example, when using Algorithm \ref{algo:l2}.I v2 for the Ijcnn1\_3500 dataset, with the value of parameter $C$ being equal to 0.1 and after a preprocessing of 28850.40 seconds, we are able to solve the problem in less than twenty seconds. While using the Algorithm of \cite{bonami} using Corollary \ref{coro:initialMl2} or the Ind. Const. (LIC) approach, the GAPs are  73.18\% and 92.10\%, respectively. We observed that we are able to solve many more problems to optimality when using  Algorithm \ref{algo:l2}.I v2 than when using Algorithm \ref{algo:l2}.II. In almost all the cases, the performance of Algorithm \ref{algo:l2}.I v2 is better than the Algorithm of \cite{bonami} and the Ind. Const. (LIC) approach. 

To summarize, we tested Algorithm \ref{algo:l2} in several datasets and in the majority of the cases it performed better than the latest algorithms published in the literature. Therefore, our strategies could be considered as useful tools for the enhancement of the \mbox{(RL-$\ell_2$-M)} formulation.

	For interested readers, a comparison between the proposed methodologies and the existing solution methods in terms of classification performance can be found in the supplementary material. 	

\section{Conclusions}
\label{conclusions}

As stated in \cite{PedroDuarte17}, the search for exact solutions of the optimization problems resulting from the ramp loss models with large datasets is an open problem. This paper has presented various new exact approaches
which are applicable to larger datasets and they have proven to be faster than the state-of-the-art algorithms.  
Unlike other resolution methods in the literature, these exact approaches do not use any auxiliary mixed integer programming model. In fact, the valid inequalites included in the formulation, the theorems that tighten the big M parameters, and the techniques for obtaining bounds on the variables only make use of auxiliary relaxed problems. Such problems allow the improvement of the relaxation of the original model. Consequentially, they enhance the behavior of the branch and bound algorithm.

There are several opportunities for future work: firstly, the presented strategies could be adapted for a $\ell_2$-norm ramp loss model with a kernel function to make the classifier more accurate. Secondly, a ramp loss model could be designed to include some constraints for feature selection. 
\section*{Acknowledgements}

Thanks to the \textit{Agencia Estatal de
	Investigaci\'on (AEI) and the European Regional Development Fund (ERDF):} project  MTM2016-74983-C2-2-R, \textit{Universidad de C\'{a}diz (UCA):}
	UCA Programa de Fomento e Impulso de la actividad Investigadora (2018) and PhD grant UCA/REC01VI/2017. Thanks to the author of \cite{31Brooks2011} for providing the datasets used in his paper and thanks to the authors of \cite{bonami} for providing the datasets and the codes of the algorithm presented in their paper. We would also like to thank the anonymous reviewers for their suggestions that have improved the quality of the paper.

\appendix
\section{Appendix}\label{sec:Appendix}

This appendix includes the proof of Theorem \ref{tm:lagrangel22}, Corollary \ref{coro:l2}, and Proposition \ref{prop:Mclusterl2+}.

\noindent\textbf{Proof of Theorem \ref{tm:lagrangel22}:}
Let $\bar{\alpha}$ be the vector of optimal values for the dual variables associated with family of constraints $\eqref{ec:linearwk}$ of problem $(\overline{\mbox{Re-RL-$\ell_2$}})_{k_0}.$ By applying Proposition 3.4.2 \cite{bertsekas1999nonlinear}, it holds that: 
\begin{align}
Z_{k_0}=& \dfrac{1}{2}\sum_{k=1}^{d}\mbox{$\bar{w}_{k}^\ast$}^{ 2} + C \left(\sum_{i=1}^{n}\bar{\xi_i}^\ast + 2\sum_{i=1}^{n} \bar{z}_i^\ast\right)+\dfrac{1}{2}\tilde{w}_{k_0}^2+\tilde{w}_{k_0} \bar{w}_{k_0}^{\ast} \nonumber\\& + \sum_{i=1}^{n} \bar{\alpha_i}\left( 1- \bar{\xi_i}^\ast -y_i\tilde{w}_{k_0} x_{ik_0}  -M_i\bar{z_i}^\ast- y_i \sum_{k=1}^{d}\bar{w}_k^{\ast} x_{ik}- y_i\bar{\bb}^\ast\right). \nonumber
\end{align}

Additionally, since $\bar{w}^{\ast}_{k_0}=0,$ we have:
\begin{align}
Z_{k_0}=& \frac{1}{2}\sum_{k=1,k\neq k_0}^{d}\mbox{$\bar{w}_{k}^\ast$}^{ 2}  + C \left(\sum_{i=1}^{n}\bar{\xi_i}^\ast + 2\sum_{i=1}^{n} \bar{z}_i^\ast\right)+\dfrac{1}{2}\tilde{w}_{k_0}^2 \nonumber\\& + \sum_{i=1}^{n} \bar{\alpha_i}\left( 1- \bar{\xi_i}^\ast -y_i\tilde{w}_{k_0} x_{ik_0}  -M_i\bar{z_i}^\ast- y_i \sum_{k=1,k\neq k_0}^{d}\bar{w}_k^{\ast} x_{ik}- y_i\bar{\bb}^\ast\right).\label{cond:lang22}
\end{align}

On the other hand, formulation $(\overline{\mbox{Re-RL-$\ell_2$}})_{k_0}$ with the additional constraint $\bar{w}_{k_0}= \hat{w}_{k_0}$ and where family of constraints \eqref{ec:linearwk} has been dualized is the following: 
\begin{alignat}{3}
& \mbox{min}&  \quad & 	\displaystyle  \frac{1}{2}\sum_{k=1}^{d}{\bar{w}_k}^2 + C  \left( \sum_{i=1}^n\bar{\xi}_i + 2 \sum_{i=1}^n \bar{z}_i \right) +\dfrac{1}{2}\tilde{w}_{k_0}^2+\tilde{w}_{k_0} \bar{w}_{k_0}  \nonumber&&\\
&&& + \sum_{i=1}^{n} \alpha_i\left( 1- \bar{\xi}_i -y_i\tilde{w}_{k_0} x_{ik_0}-M_i\bar{z}_i- y_i \sum_{k=1}^{d}\bar{w}_k x_{ik}- y_i\bar{\bb}\right),   \nonumber&& \\
&\mbox{s.t.} &  \quad & \eqref{ec:cond1prime} -\eqref{ec:zprime},  && \nonumber \\
& & & \bar{w}_{k_0}= \hat{w}_{k_0},\nonumber   && 
\end{alignat}
where $\alpha_i\geq 0.$ Therefore, this problem can be rewritten as follows, 
\begin{alignat}{4}
(\overline{\mbox{Lg-RL-$\ell_2$}})_{k_0} &\quad &\mbox{min} \quad &   	\displaystyle  \frac{1}{2} \sum_{k=1,k\neq k_0}^{d}{\bar{w}_k}^2 + C  \left( \sum_{i=1}^n\bar{\xi}_i + 2 \sum_{i=1}^n \bar{z}_i\right)+ \hat{w}_{k_0}\left(\dfrac{1}{2}\hat{w}_{k_0}+\tilde{w}_{k_0} - \sum_{i=1}^{n} \alpha_iy_ix_{ik_0} \right)   &&\nonumber\\ &&& \dfrac{1}{2}\tilde{w}_{k_0}^2+\sum_{i=1}^{n} \alpha_i\left( 1- \bar{\xi}_i -y_i\tilde{w}_{k_0} x_{ik_0} -M_i\bar{z}_i- y_i \sum_{k=1, k\neq k_0}^{d}\bar{w}_k x_{ik}- y_i\bar{\bb}\right) , & & \nonumber \\
& & \mbox{s.t.} \quad& \eqref{ec:cond1prime}-\eqref{ec:zprime}.\nonumber 
\end{alignat}

Note that ($\bar{w}^{\ast}, \bar{\bb}^\ast, \bar{\xi}^\ast, \bar{z}^\ast$), an optimal solution of $(\overline{\mbox{Re-RL-$\ell_2$}})_{k_0}$, is feasible for the problem above. Additionally, any feasible solution of problem $(\overline{\mbox{Lg-RL-$\ell_2$}})_{k_0},$ when taking $\bar{w}_{k_0}=0,$ is feasible for $(\overline{\mbox{Re-RL-$\ell_2$}})_{k_0}$ where family of constraints \eqref{ec:linearwk}  has been dualized. Therefore, for $\alpha= \bar{\alpha}$ and when using \eqref{cond:lang22}, the optimal objective value of the problem above is $Z_{k_0}+ \hat{w}_{k_0}\left(\dfrac{1}{2}\hat{w}_{k_0}+\tilde{w}_{k_0} - \displaystyle \sum_{i=1}^{n} \bar{\alpha_i}y_ix_{ik_0} \right).$ This is a lower bound of the optimal value of $(\overline{\mbox{Re-RL-$\ell_2$}})_{k_0}$ with the additional constraint that $\bar{w}_{k_0}= \hat{w}_{k_0}.$

\fin

\noindent\textbf{Proof of Corollary \ref{coro:l2}:}
An equivalent model $(\overline{\mbox{Re-RL-$\ell_2$}})_{k_0}$ can be built from an optimal solution of (Re-RL-$\ell_2$), $({w}^*, \bb^*, \xi^*, z^*),$ which satisfies that ${w_{k_0}}^*=\tilde{w}_{k_0}$. In this situation, both models have the same optimal objective value. Since model (Re-RL-$\ell_2$) is the relaxation of model (RL-$\ell_2$-M), an upper bound of (RL-$\ell_2$-M) would be an upper bound of (Re-RL-$\ell_2$).
Thus, if the objective value of an optimal solution of  $(\overline{\mbox{Re-RL-$\ell_2$}})_{k_0}$ restricting $\bar{w}_{k_0}= \hat{w}_{k_0}$ is bigger than $\mbox{UB}_{\tiny{\mbox{RL-}\ell_2}},$ the value $\bar{w}_{k_0}= \hat{w}_{k_0}$ can be discarded as optimal solution of (RL-$\ell_2$-M). This is because any solution with this value will provide a solution whose objective value is worse than $\mbox{UB}_{\tiny{\mbox{RL-}\ell_2}}$. Therefore, we can restrict ourselves to the values of $\hat{w}_{k_0}$ in such a way that $Z_{\hat{k}_0} \leq \mbox{UB}_{\tiny{\mbox{RL-}\ell_2}}$. Therefore, according to Theorem \ref{tm:lagrangel22}, it holds that: 
$$Z_{k_0}+ \bar{w}_{k_0}\left(\dfrac{1}{2}\bar{w}_{k_0}+\tilde{w}_{k_0} - \displaystyle\sum_{i=1}^{n} \bar{\alpha_i}y_ix_{ik_0} \right)\leq \mbox{UB}_{\tiny{\mbox{RL-}\ell_2}}.$$

Therefore, $\bar{w}_{k_0}$ verifies:
\begin{equation} \nonumber
\displaystyle \sum_{i=1}^{n} \bar{\alpha_i}y_ix_{ik_0}-\tilde{w}_{k_0} - \sqrt{\displaystyle \left(\tilde{w}_{k_0}-\sum_{i=1}^{n} \bar{\alpha_i}y_ix_{ik_0}  \right)^2 -2 \left( Z_{k_0}-\mbox{UB}_{\tiny{\mbox{RL-}\ell_2}}\right)}\leq \bar{w}_{k_0},
\end{equation}
\begin{equation}  \nonumber
\bar{w}_{k_0} \leq  \displaystyle \sum_{i=1}^{n} \bar{\alpha_i}y_ix_{ik_0}-\tilde{w}_{k_0} + \sqrt{\displaystyle \left(\tilde{w}_{k_0}-\sum_{i=1}^{n} \bar{\alpha_i}y_ix_{ik_0}  \right)^2 -2 \left( Z_{k_0}-\mbox{UB}_{\tiny{\mbox{RL-}\ell_2}}\right)}.
\end{equation}
Taking into account that $\bar{w}_{k_0}=w_{k_0}-\tilde{w}_{k_0},$ the result is obtained. \fin

\noindent\textbf{Proof of Proposition \ref{prop:Mclusterl2+}:}
For each $k\in D$ and $i \in c_+$ when $y_{i}=+1,$ the following inequalities are satisfied:  $-w_kx_{i_0k} \leq |w_k||\bar{x}_{c+k}|$. Taking the summation in $k,$ we have: 
\begin{equation}1 - \left(\sum_{k=1}^{d}w_k x_{ik}+\bb \right)\leq 1+ \left(\sum_{k=1}^{d} |w_k||\bar{x}_{c+k}| -\bb \right),\; \; \mbox{for } i\in c_+. \end{equation}

Since $0\leq \xi_{i}\leq2$, for $i \in c_+$ and the $v$-variables model the absolute value of the $w$-variables, we can conclude:
\begin{equation}1 -\xi_{i}- \left(\sum_{k=1}^{d}w_k x_{i_{0}k}+\bb \right)\leq 1+ \left(\sum_{k=1}^{d} v_k|\bar{x}_{c+k}| -\bb \right),\; \; \mbox{for } i\in c_+. \nonumber\end{equation}

 Therefore, the optimal solution of problem $\left(\mbox{UB}^{\; \ell_2}_{M_{c_+}}\right)$ is an upper bound of problem $\left(\mbox{UB}^{\;\ell_2}_{M_{i_0}}\right)$, for $i\in c_+$. Similarly, it can be proven that the optimal solution of problem $\left(\mbox{UB}^{\; \ell_2}_{M_{c_-}}\right)$ is an upper bound of problem $\left(\mbox{UB}^{\;\ell_2}_{M_{i_0}}\right)$, for $i\in c_-$. Therefore, according to Proposition \ref{propDirectl2}, the result holds.
 
\fin

\end{document}